\documentclass[mathpazo]{cicp}
\usepackage{graphicx, amsmath, amssymb}
\usepackage[american]{babel}
\usepackage{epsfig}
\usepackage{psfrag}
\usepackage{makeidx}

\newtheorem{thm}{Theorem}[section]

\begin{document}
\def\bnu{\mbox{\boldmath{$\nu$}}}
\def\bphi{\mbox{\boldmath{$\phi$}}}
\def\btheta{\mbox{\boldmath{$\vartheta$}}}
\def\bphiH{\mbox{\boldmath{$\hat{\phi}$}}}
\def\bthetaH{\mbox{\boldmath{$\hat{\vartheta}$}}}
\def\bPhi{\mbox{\boldmath{$\Phi$}}}
\def\bphiT{\mbox{\underline{\underline{\boldmath{$\phi$}}}} }
\def\bpsi{\mbox{\boldmath{$\psi$}}}
\def\bPsi{\mbox{\boldmath{$\Psi$}}}
\def\bg{\mbox{\boldmath{$g$}}}
\def\bs{\mbox{\boldmath{$s$}}}
\def\bu{\mbox{\boldmath{$u$}}}
\def\bw{\mbox{\boldmath{$w$}}}
\def\bbf{\mbox{\boldmath{$f$}}}
\def\btu{\mbox{\boldmath{$\tilde{u}$}}}
\def\bz{\mbox{\boldmath{$z$}}}
\def\bM{\mbox{\boldmath{$\mathcal{M}$}}}
\def\bA{\mbox{\boldmath{$\mathcal{A}$}}}
\def\bF{\mbox{\boldmath{$\mathcal{F}$}}}
\def\btF{\mbox{\boldmath{$\tilde{\mathcal{F}}$}}}
\def\bG{\mbox{\boldmath{$\mathcal{G}$}}}
\def\bK{\mbox{\boldmath{$\mathcal{K}$}}}
\def\vdiv{\mathbf{div}\,}
\def\vcurl{\mbox{$\mathbf{curl}$\,}}
\def\curl{\mbox{curl\,}}
\def\div{\mbox{div\,}}
\def\det{\mbox{det}}
\def\vgrad{\mathbf{grad\,}}
\def\qT{\mbox{\underline{\underline{\boldmath{$q$}}}}}
\def\bv{\mbox{\boldmath{$v$}}}
\def\bx{\mbox{\boldmath{$x$}}}
\def\bk{\mbox{$\mathbf{k}$}}
\def\bs{\mbox{$\mathbf{s}$}}
\def\velocity{\mbox{$\mathbf{v}$}}
\def\Hz{\mbox{Hz\,}}
\def\m{\mbox{m}}
\def\e{\mbox{e}}
\def\a{\mbox{a}}
\def\s{\mbox{s}}
\def\V{\mbox{\boldmath{$V$}}}
\def\U{\mbox{\boldmath{$U$}}}
\def\bw{\mbox{\boldmath{$w$}}}
\title{Efficient PML for the wave equation}
\author[Grote M J et.~al.]{Marcus J. Grote\affil{1}, and  Imbo Sim\affil{2}\corrauth}
\address{\affilnum{1}\ Department of Mathematics, University of Basel, CH-4051 Basel, Switzerland \\
\affilnum{2}\ Institut d'Analyse et Calcul Scientifique, Ecole Polytechnique F\'ed\'erale de Lausanne, CH-1015 Lausanne, Switzerland}
\emails{{\tt Marcus.Grote@unibas.ch}(Marcus~J.~Grote), {\tt Imbo.Sim@epfl.ch} (Imbo Sim)}
\begin{abstract}
In the last decade, the perfectly matched layer (PML) approach 
 has proved a flexible and accurate method for the simulation of waves in unbounded media. 
 Most PML formulations, however, usually require wave equations stated in
 their standard second-order form to be reformulated as first-order systems, thereby
introducing many additional unknowns. 
To circumvent this cumbersome and somewhat expensive step, we instead propose
a simple PML formulation directly for the wave
 equation in its second-order form. Inside the absorbing layer,
our formulation requires only two auxiliary variables in two space dimensions and
four auxiliary variables in three space dimensions; hence it is cheap to implement.
Since our formulation requires no higher derivatives, it is also easily coupled 
with standard finite difference or finite element methods. Strong stability is
proved while numerical examples in two and three space dimensions
illustrate the accuracy and long time stability of our PML formulation.
\end{abstract}
\maketitle
\section{Introduction}
The accurate and reliable simulation of wave propagations in unbounded media is of fundamental importance in a wide range of applications.
The perfectly matched layer (PML) approach \cite{ber} has proved a flexible and accurate method for the simulation of waves. It consists in surrounding the computational domain by an absorbing layer, which 
generates no reflections at its interface with the computational domain; hence, it is perfectly matched.
Inside the absorbing layer a damping term is added to the wave equation, which acts only in 
the direction perpendicular to the layer.
This approach is analogous to the physical treatment of the walls of an anechoic chamber and provides an alternative to absorbing or nonrelfecting boundary conditions \cite{GrKe1, GrKe2,Ha,hh98,hmg}. \\
The initial PML formulation of B\'erenger \cite{ber} was based on splitting the electromagnetic fields into two parts, 
the first containing the tangential derivatives and the second containing the normal derivatives;
damping was then enforced only upon the normal component. 
Later Abarbanel and Gottlieb \cite{agh} showed that B\'erenger's approach 
was only weakly well-posed due to the unphysical splitting of the field variables. 
Several strongly well-posed approaches have been suggested since, 
some of which were shown to be linearly equivalent \cite{ahk, zc}. \\
 The PML approach has proved very successful in practice, because of its simplicity, versatility, and
robust treatment of corners. Once discretized and truncated at a finite thickness, the layer is no
longer perfectly absorbing and the optimal damping parameters need to be determined via numerical experiments.
Stability properties of the PML approach has been analyzed in several works, such as in \cite{dj,agh,ahk,bfj}
among others.\\
The best implementation in the time domain is still under debate. Most PML formulations 
require wave equations stated in their standard second-order form to be reformulated as first-order hyperbolic 
systems, thereby introducing many additional unknowns. 
Here we propose instead a simple PML formulation directly for the second-order wave equation 
both in two and in three space dimensions. Our formulation also requires fewer auxiliary variables
than previous formulations for the second-order wave equation -- see \cite{ak2,c,sp}, for instance.\\
Our paper is organized as follows. In Section 2 we derive a PML formulation for the
wave equation in its standard second-order form. By judiciously choosing the auxiliary variables
in the Laplace transformed domain, the resulting PML modified equations require only
two auxiliary variables in two dimensions and four auxiliary variables in three dimensions
inside the absorbing layer. Next, in Section 3 we prove stability of our PML formulation
by using standard theory from \cite{KL}. The finite difference discretization
of the PML modified wave equation is shown in Section 4. In Section 5, our numerical results 
both in two and three space dimensions demonstrate the accuracy and long time stability of the PML formulation.
\section{PML formulation}
We consider a time dependent wave field $u$ propagating through unbounded three dimensional space and assume that all sources and initial disturbances are confined to the rectangular domain $\Omega = [-a_1, a_1]\times[-a_2, a_2]\times[-a_3, a_3]$,  $ a_1, \, a_2, \, a_3 > 0$. Outside $\Omega$, we further assume the speed of propagation $c > 0$ to be constant; hence, all waves are purely outgoing in the unbounded exterior $\mathbb{R}^3 \backslash \Omega$. Inside $\Omega$, the wave field $u(x_1,x_2,x_3,t)$ satisfies
\begin{alignat}2
\label{eq:wave1}
    u_{tt} - \nabla\cdot \left(\, c^2 \, \nabla u\right) &= f
    & \qquad &  t > 0,\\[0.1cm]
    \label{eq:wave3}
    u & = u_0 & \qquad & t = 0 ,\\[0.1cm]
    \label{eq:wave4}
     u_t  &= v_0 & \quad & t = 0 .
\end{alignat}
We wish to truncate the unbounded exterior and thereby restrict the computation to the
finite computational domain $\Omega$. In doing so, we need to ensure that all waves propagating
outward leave $\Omega$ without spurious reflection. Thus we shall surround $\Omega$ by a perfectly
matched layer (PML) of thickness $L_i$, $i = 1,2,3$, in each coordinate which is designed
to absorb the waves exiting $\Omega$. Inside the absorbing layer, $u$ then satisfies a 
modified wave equation whose solutions decay exponentially fast with distance from 
the computational domain. \\
Following \cite{agh, ahk}, we let $\hat{u}$ denote the Laplace transform of $u$, defined
as 
\begin{equation}
\label{eq:laplace_trafo}
\hat{u} = \hat{u} \, (\bx, \s) = \int_0^{\infty} \e^{\mathrm{s} \, t} \, u(\bx,t) \, dt, 
\qquad  s \in \mathbb{C} .
\end{equation}
Outside $\Omega$, $\hat{u}$ then satisfies the Helmholtz equation,
\begin{equation}
\label{eq:helmholtz1}
\s^2 \, \hat{u} =  \frac{\partial}{\partial x_1} \bigg( c^2 \, \frac{\partial \hat{u}}{\partial x_1} \bigg)
+ \frac{\partial}{\partial x_2} \bigg( c^2 \, \frac{\partial \hat{u}}{\partial x_2}\bigg)
+ \frac{\partial}{\partial x_3} \bigg( c^2 \, \frac{\partial \hat{u}}{\partial x_3}\bigg) .
\end{equation}
Next, we introduce the coordinate transformation  
\begin{equation}
\label{eq:traf}
x_i   \mapsto \tilde{x}_i:= x_i + \frac{1}{\s} \int_0^{x_i} \zeta_i(x) \, d x,  \qquad   i=1,2,3 ,
\end{equation}
where the damping profile $\zeta_i$ is positive inside the absorbing layer, 
$|x_i| > a_i$, $i = 1,2,3$, but vanishes inside $\Omega$. If we now require $\hat{u}$ 
to satisfy the modified Helmholtz equation in those stretched coordinates, 
\begin{equation}
\label{eq:helmholtz1_modification}
s^2 \, \hat{u} =
\frac{\partial}{\partial
  \tilde{x}_1} \bigg( c^2 \, \frac{\partial \hat{u}}{\partial \tilde{x}_1} \bigg)
+ \frac{\partial}{\partial \tilde{x}_2} \bigg( c^2 \, \frac{\partial \hat{u}}{\partial \tilde{x}_2} \bigg)
+ \frac{\partial}{\partial \tilde{x}_3} \bigg( c^2 \, \frac{\partial \hat{u}}{\partial \tilde{x}_3} \bigg),
\end{equation}
it is well-known that $u$ will remain unaltered inside $\Omega$, but decay exponentially fast inside the layer; 
hence the absorbing layer will be perfectly matched. In fact, the (unsplit) PML modified Helmholtz equation
(\ref{eq:helmholtz1_modification}) in the Laplace transformed domain is standard \cite{agh, ahk}. 
The difficulty lies in transforming
(\ref{eq:helmholtz1_modification}) back to the time domain, without introducing 
high order derivatives or too many auxiliary variables.  \\
From $(\ref{eq:traf})$, we observe that partial differentiation with respect to $\tilde{x}_i$ is related
to partial differentiation with respect to the physical coordinate, $x_i$, through
\begin{equation}
\label{eq:tsfm}
\frac{\partial}{\partial \tilde{x}_i} = \frac{\s}{\s + \zeta_i} \frac{\partial}{\partial x_i}.
\end{equation}
We now let $\gamma_i =\gamma_i(\zeta_i ; \, \,  \s)$, $i=1,2,3$ denote
\begin{equation}
\label{eq:gamma}
\gamma_i := 1 + \frac{\zeta_i}{\s} .
\end{equation}
Then, by replacing partial derivatives according to (\ref{eq:tsfm}) and multiplying the resulting 
expression by $\gamma_1 \, \gamma_2 \, \gamma_3 $, we rewrite (\ref{eq:helmholtz1_modification})  in physical
coordinates as 
\begin{equation}
\label{eq:helmholtz2}
\s^2 \, \gamma_1 \,  \gamma_2 \,  \gamma_3\, \hat{u} =  \frac{\partial}{\partial x_1} \bigg(c^2 \, \frac{\gamma_2 \, \gamma_3 }{\gamma_1}  \, \frac{\partial \hat{u}}{\partial x_1} \bigg) + \frac{\partial}{\partial x_2} \bigg(c^2  \, 
\frac{\gamma_3 \, \gamma_1 }{\gamma_2} \, \frac{\partial \hat{u}}{\partial x_2} \bigg) + \frac{\partial}{\partial x_3} \bigg(c^2  \, \frac{\gamma_1 \, \gamma_2 }{\gamma_3} \, \frac{\partial \hat{u}}{\partial x_3} \bigg)          .
\end{equation}
From (\ref{eq:gamma}) we derive after some algebra the following identities:
\begin{equation}
\label{eq:gamma_algbra}
\begin{split}
\frac{\gamma_2 \, \gamma_3 }{\s \, \gamma_1} &= 1 + \frac{(\zeta_2 + \zeta_3 - \zeta_1 )\s + \zeta_2 \zeta_3 }{(\s + \zeta_1)\s},\\
\frac{\gamma_3 \, \gamma_1 }{\s \, \gamma_2} &= 1 + \frac{(\zeta_3 + \zeta_1 - \zeta_2 )\s + \zeta_3 \zeta_1 }{(\s + \zeta_2)\s},\\
\frac{\gamma_1 \, \gamma_2 }{\s \, \gamma_3} &= 1 + \frac{(\zeta_1 + \zeta_2 - \zeta_3 )\s + \zeta_1 \zeta_2 }{(\s + \zeta_3)\s}.
\end{split}
\end{equation}
By using (\ref{eq:gamma_algbra}) in (\ref{eq:helmholtz2}) we find
\begin{equation}
\label{eq:hholtz3}
\begin{split}
\quad &\, \left(\s^2 + \s \, (\zeta_1 + \zeta_2 + \zeta_3) + (\zeta_1 \, \zeta_2 + \zeta_2 \, \zeta_3 + \zeta_3 \, \zeta_1 )
+ \frac{\zeta_1 \, \zeta_2 \, \zeta_3 }{\s} \right) \, \hat{u} \\
\quad &\hspace{-0.15cm}= \frac{\partial}{\partial x_1} \bigg( c^2 \, \frac{\partial \hat{u}}{\partial x_1} \bigg)
+ \frac{\partial}{\partial x_2} \bigg( c^2 \, \frac{\partial \hat{u}}{\partial x_2} \bigg)
+ \frac{\partial}{\partial x_3} \bigg( c^2 \, \frac{\partial \hat{u}}{\partial x_3} \bigg) \\
\quad &\hspace{-0.15cm}
+\frac{\partial}{\partial x_1} \bigg( c^2\, \Big( \frac{(\zeta_2 + \zeta_3 - \zeta_1)\s + \zeta_2 \, \zeta_3}{(\s +\zeta_1)\s} \Big)  \frac{\partial \hat{u}}{\partial x_1} \bigg) 
+\frac{\partial}{\partial x_2} \bigg( c^2\, \Big( \frac{(\zeta_3 + \zeta_1 - \zeta_2)\s + \zeta_3 \, \zeta_1}{(\s +\zeta_2)\s} \Big)  \frac{\partial \hat{u}}{\partial x_2} \bigg) \\
\quad &\hspace{-0.15cm}
+\frac{\partial}{\partial x_3} \bigg( c^2\, \Big( \frac{(\zeta_1 + \zeta_2 - \zeta_3)\s + \zeta_1 \, \zeta_2}{(\s +\zeta_3)\s} \Big)  \frac{\partial \hat{u}}{\partial x_3} \bigg) .
\end{split}
\end{equation}
Next, we introduce the auxiliary functions $\psi$ and $\bphi = \left(\phi_1, \phi_2, \phi_3 \right)^{\top} $,
\begin{equation*}
\begin{split}
\widehat{\psi} &= \frac{1}{\s} \hat{u}, \\
\widehat{\phi}_1 &= c^2 \,\bigg(  \frac{\zeta_2 + \zeta_3 - \zeta_1}{\s+ \zeta_1} + \frac{\zeta_2 \, \zeta_3 }{(\s+ \zeta_1)\s}  \bigg) \frac{\partial \hat{u}}{\partial x_1}, \\
\widehat{\phi}_2 &= c^2 \,\bigg(  \frac{\zeta_3 + \zeta_1 - \zeta_2}{\s+ \zeta_2} + \frac{\zeta_3 \, \zeta_1 }{(\s+ \zeta_2)\s}  \bigg)   \frac{\partial \hat{u}}{\partial x_2}, \\ 
\widehat{\phi}_3 &= c^2 \,\bigg(  \frac{\zeta_1 + \zeta_2 - \zeta_3}{\s+ \zeta_3} + \frac{\zeta_1 \, \zeta_2 }{(\s+ \zeta_3)\s}  \bigg)   \frac{\partial \hat{u}}{\partial x_3},
\end{split}
\end{equation*}
or equivalently
\begin{equation*}
\begin{split}
\s \, \widehat{\psi} &= \hat{u}, \\
(\s + \zeta_1)\, \widehat{\phi}_1 &= c^2 \bigg( \left( \zeta_2 + \zeta_3 - \zeta_1 \right) + \frac{\zeta_2 \, \zeta_3}{\s}     \bigg) \frac{\partial \hat{u}}{\partial x_1}, \\
(\s + \zeta_2)\, \widehat{\phi}_2 &= c^2 \bigg( \left( \zeta_3 + \zeta_1 - \zeta_2 \right) + \frac{\zeta_3 \, \zeta_1}{\s}     \bigg) \frac{\partial \hat{u}}{\partial x_2}, \\
\text{and} \quad (\s + \zeta_3)\, \widehat{\phi}_3 &= c^2 \bigg( \left( \zeta_1 + \zeta_2 - \zeta_3 \right) + \frac{\zeta_1 \, \zeta_2}{\s}     \bigg) \frac{\partial \hat{u}}{\partial x_3}.
\end{split}
\end{equation*}
Finally, we use the above relations in (\ref{eq:hholtz3}) and transform the resulting equations back
to the time domain, which yields the PML modified wave equation
\begin{equation}
\label{pml-eq_3d}
\begin{split}
u_{tt}
+ \left(\zeta_1 + \zeta_2 + \zeta_3 \right)  u_t
+ \left(\zeta_1 \, \zeta_2  +  \zeta_2 \, \zeta_3  + \zeta_3 \, \zeta_1\right)u 
& = \nabla\cdot \left(\, c^2 \, \nabla u \right) +  \nabla \cdot \bphi - \zeta_1 \, \zeta_2 \, \zeta_3 \, \psi, \\
\bphi_t &= \Gamma_1 \, \bphi + c^2 \, \Gamma_2 \, \nabla u + c^2 \, \Gamma_3 \, \nabla \psi , \\
\psi_t  &= u,
\end{split}
\end{equation}
where
\begin{equation*}
\begin{split}
\Gamma_1 &= \begin{bmatrix}
-\zeta_1 &    0     & 0 \\
    0    & -\zeta_2 & 0 \\
    0    &    0     & -\zeta_3 
\end{bmatrix}, \quad
\Gamma_2 = \begin{bmatrix}
\zeta_2 + \zeta_3 - \zeta_1 & 0 & 0 \\
0 & \zeta_3 + \zeta_1 - \zeta_2 & 0 \\
0 & 0 & \zeta_1 + \zeta_2 - \zeta_3
\end{bmatrix} \\
\text{and}\quad \Gamma_3 &= \begin{bmatrix}
\zeta_2 \, \zeta_3 & 0 & 0 \\
 0 & \zeta_3 \,  \zeta_1 & 0 \\
 0 & 0 & \zeta_1 \, \zeta_2
\end{bmatrix} .
\end{split}
\end{equation*}
In the interior of $\Omega$, the damping profiles $\zeta_i, \, i=1,2,3$ and the auxiliary variables $\bphi$, $\psi$
vanish; hence, (\ref{pml-eq_3d})
reduces to (\ref{eq:wave1}) in $\Omega$. 
Because our PML formulation (\ref{pml-eq_3d}) requires only four auxiliary scalar variables 
$\phi_1, \, \phi_2, \, \phi_3,  \, \psi$ inside the layer and no high order derivatives,
its implementation is not only straightforward but also cheap to implement. \\
In two space dimensions, $\zeta_3$ and $\phi_3$ and $\psi$ vanish and our PML formulation reduces to
\begin{equation}
\label{pml-eq_2d}
\begin{split}
u_{tt} + \left(\zeta_1 + \zeta_2 \right)  u_t + \zeta_1 \, \zeta_2  u & = \nabla\cdot \left(\, c^2 \, \nabla u \right) +  \nabla \cdot \bphi ,\\
\bphi_t &= \Gamma_1 \, \, \bphi +  c^2 \, \Gamma_2 \, \nabla u , 
\end{split}
\end{equation}
where
\begin{equation*}
\Gamma_1 = \begin{bmatrix}
- \zeta_1 & 0 \\
0 &  - \zeta_2
\end{bmatrix}, \quad
\Gamma_2 = \begin{bmatrix}
\zeta_2 - \zeta_1 & 0 \\
0 & \zeta_1 - \zeta_2
\end{bmatrix} .
\end{equation*}
Remarkably only two auxiliary functions are needed here. \\
The choice of the damping profiles $\zeta_i(x) \geq 0,\quad i=1,2,3$ is arbitrary; 
it can be constant, linear, or quadratic among others. In our computations, we always use
\begin{equation}
\label{eq:zeta}
\zeta_i(x_i) = 
\begin{cases} 0 & \text{for } |x_i| < a_i, \quad i=1,2,3   \\ 
\bar{\zeta}_i \left( \frac{|x_i-a_i|}{L_i} - \frac{\sin \left( \frac{2 \pi \, |x_i-a_i|}{L_i} \right)}{2 \pi} \right) 
& \text{for }  a_i \leq |x_i| \leq a_i + L_i, \quad i=1,2,3.
\end{cases}
\end{equation}
Because $\zeta_i(x)$ is twice continuously differentiable throughout the interface at $|x_i| = a_i$, no
special transmission conditions are needed there. The constant $\bar{\zeta}_i$ depends on the 
discretization and the thickness of the layer, which in practice is truncated by a homogeneous 
Dirichlet (or Neumann) boundary condition. Then the relative reflection, $R$, is given by
\begin{equation}
\bar{\zeta}_i = \frac{c}{L_i} \log \Big(\frac{1}{R} \Big), \quad i=1,2,3 .
\end{equation}
In Figure 1 we show damping profiles for different values of $\bar{\zeta}_i$.
\begin{figure}[htbp]
\begin{center}
\epsfig{file = 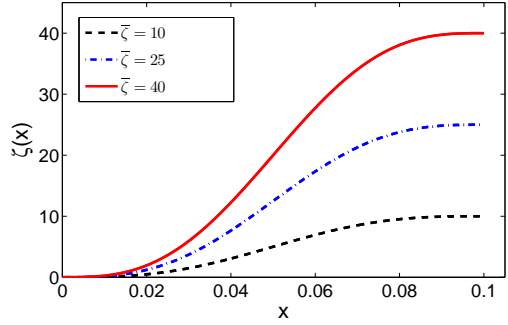, width=7.5cm}  
\end{center}
\label{fig:damp}
\caption{The damping profile $\zeta_i(x_i)$ given by (\ref{eq:zeta}) is shown for different values of $\bar{\zeta}_i$,
with $c=1$ and $L_i = 0.1$~.}
\end{figure}
\section{Stability}
We now establish the stability and well-posedness of our PML formulation, first in two and
then in three space dimensions, where we assume that the absorbing layer extends to infinity.
Here we follow standard stability theory for hyperbolic systems \cite{KL}, which we briefly recall below.  \\
Consider a general Cauchy problem,
\begin{equation}
U_t = P \left(\frac{\partial }{\partial x} \right) U,   \quad  0 \leq t \leq T, \quad U \in \mathbb{R}^p,
\end{equation}
where $P(\partial_x)$ denotes a linear differential operator,
with initial conditions
\begin{equation}
U(x,0) = U_0(x), \quad x \in \mathbb{R}^3.
\end{equation}
Following \cite{KL}, the Cauchy Problem 
is \textit{weakly (resp. strongly) well-posed}, if the solution $U( \, \cdot  \,, \, t)$ satisfies
\begin{equation}
\| U( \, \cdot \, , \, t) \|_{L_2} \leq K e^{\alpha t} \|U( \, \cdot \, , \, 0) \|_{H^s} 
\end{equation}
with $s>0$ (resp. $s = 0$). 
The Cauchy Problem is \textit{weakly (resp. strongly) stable}, if the solution $U( \, \cdot  \,, \, t)$ satisfies
\begin{equation}
\| U( \, \cdot \, , \, t) \|_{L_2} \leq K \, (1+t)^s  \|U( \, \cdot \, , \, 0) \|_{H^s} 
\end{equation}
with $s>0$ (resp. $s = 0$).
A necessary and sufficient condition for weak \textit{well-posedness (resp. stability)} 
is that all eigenvalues $\lambda$ of the operator $P\left(i k)\right)$ satisfy
\begin{equation}
\Re\{\lambda \left(P(i k)\right)\} \leq C, \qquad k\in\mathbb{R},
\end{equation}
with $C >0$ (resp. $C = 0$) independent of $k$.
For strong \textit{well-posedness (resp. stability)}, the corresponding eigenvectors must also be
complete. \\
By rewriting the PML-modified wave equations (\ref{pml-eq_3d}), (\ref{pml-eq_2d}) 
as a first-order hyperbolic system and applying the stability theory from \cite{KL} delineated
above, we can prove the following two stability results.
{\thm{ The Cauchy problem for the PML formulation (\ref{pml-eq_2d}) in two space dimensions
is strongly stable for $\zeta_1, \, \zeta_2 \geq 0$.}} \\
{\textbf{proof}})
\newline
For simplicity, we assume that $\zeta_1, \, \zeta_2$ are constant; note, however,
that the stability theory from \cite{KL} extends to smoothly varying coefficients.
We introduce the new variable $\mathbf{v}$ to rewrite the first equation in (\ref{pml-eq_2d})
equivalently as
\begin{equation}
\label{ref_sys1}
u_t =  -\zeta_2 \, u + \mathrm{div}\, \mathbf{v}, \quad    \mathbf{v}_t  =  -\zeta_1 \, \mathbf{v}  + c^2 \, \nabla u + \bphi \, .
\end{equation}
By using (\ref{ref_sys1}), we now rewrite (\ref{pml-eq_2d}) as a first order hyperbolic system:
\begin{equation}
\label{pml_sys_2d}
U_t = A \,U_x +  B \, U_y +  C,
\end{equation}
where 
\begin{equation}
U_t = (u, \, \phi_1,\,  \phi_2,\,  \mathbf{v}_1,\,  \mathbf{v}_2 )^{\top} , 
\end{equation}
\begin{equation}
\begin{split}
A &= \begin{bmatrix}
  0  & 0 &  0  & 1  & 0  \\
 c^2 \, (\zeta_2-\zeta_1)  & 0 &  0  & 0  & 0  \\
  0  & 0 &  0  & 0  & 0  \\      
 c^2   & 0 &  0  & 0  & 0  \\ 
  0  & 0 &  0  & 0  & 0  
\end{bmatrix}, \quad
B = \begin{bmatrix}
  0  & 0 &  0  & 0  & 1  \\
 c^2 \, (\zeta_1-\zeta_2)  & 0 &  0  & 0  & 0  \\
  0  & 0 &  0  & 0  & 0  \\      
  0   & 0 &  0  & 0  & 0  \\ 
 c^2  & 0 &  0  & 0  & 0  
\end{bmatrix}, \\ 
\text{and} \quad  
C &= - \begin{bmatrix}
 \zeta_2  &   0     &  0     & 0       & 0  \\
     0    & \zeta_1 &  0     & 0       & 0  \\
     0    &   0     &\zeta_2 & 0       & 0  \\      
     0    &   0     &  0     & \zeta_1 & 0  \\ 
     0    &   0     &  0     & 0       & \zeta_1   
\end{bmatrix}.
\end{split}
\end{equation}
By using a symbolic algebra program we find that the eigenvalues of 
the principal part of $P(i k)$ for (\ref{pml_sys_2d}) are 
\begin{equation}
\lambda \left( P \left(i k \right) \right) = \pm \, i \, c \, (k_1^2 + k_2^2)^{1\slash 2}.
\end{equation}
Thus,
\begin{equation}
\Re \{\lambda \left( P \left(i  k \right) \right) \} = 0, 
\end{equation}
while the corresponding eigenvectors are also complete for all $\zeta_1, \, \zeta_2 \geq 0$. Therefore, since $C$ is a diagonal matrix with negative entries for $\zeta_1, \, \zeta_2 \geq 0$, 
we conclude that (\ref{pml-eq_2d}) is strongly stable. \\
{\thm 
The Cauchy problem for the PML formulation (\ref{pml-eq_3d}) in three space dimensions
is strongly stable, if at least two $\zeta_j = 0, \, j = 1, 2,3$, and  weakly stable, otherwise.} \\
{\textbf{proof}})
\newline
We introduce the new variable $\mathbf{v}$ to rewrite the first equation in (\ref{pml-eq_3d})
as 
\begin{equation}
\label{ref_sys2}
u_t =  -\zeta_2 u   +  \mathrm{div}\, \mathbf{v} - \zeta_3 \psi, \quad    
\mathbf{v}_t  =  -\zeta_1 \, \mathbf{v} +  c^2 \, \nabla u + \bphi \, .
\end{equation}
By using (\ref{ref_sys2}), we can rewrite (\ref{pml-eq_3d}) as a first order hyperbolic system:
\begin{equation}
\label{pml_sys_3d}
U_t = A \,U_x +  B \, U_y +  C \, U_z + D,
\end{equation}
where 
\begin{equation}
U_t = (u, \, \phi_1,\,  \phi_2,\,  \phi_3, \,  \mathbf{v}_1,\,  \mathbf{v}_2, \,  \mathbf{v}_3, \, \psi)^{\top} , 
\end{equation}
\begin{equation}
A= \begin{bmatrix}
 0   & 0 &  0 &  0 & 1  & 0 &  0 &  0 \\
 c^2 \, (\zeta_2 + \zeta_3 - \zeta_1)  & 0 &  0  & 0 & 0 & 0 & 0 &   \zeta_2 \, \zeta_3 \\
 0   & 0 &  0 &  0 & 0  & 0 &  0 &  0 \\
 0   & 0 &  0 &  0 & 0  & 0 &  0 &  0 \\
 c^2 & 0 &  0 &  0 & 0  & 0 &  0 &  0 \\
 0   & 0 &  0 &  0 & 0  & 0 &  0 &  0 \\
 0   & 0 &  0 &  0 & 0  & 0 &  0 &  0 \\
 0   & 0 &  0 &  0 & 0  & 0 &  0 &  0 
\end{bmatrix},
\end{equation}
\begin{equation}
B = \begin{bmatrix}
 0   & 0 &  0 &  0 & 0  & 1 &  0 &  0 \\
 0   & 0 &  0 &  0 & 0  & 0 &  0 &  0 \\
 c^2 \, (\zeta_3 + \zeta_1 - \zeta_2)  & 0 &  0  & 0 & 0 & 0 &  0 &  \zeta_3 \, \zeta_1 \\
 0   & 0 &  0 &  0 & 0  & 0 &  0 &  0 \\
 0   & 0 &  0 &  0 & 0  & 0 &  0 &  0 \\
 c^2 & 0 &  0 &  0 & 0  & 0 &  0 &  0 \\
 0   & 0 &  0 &  0 & 0  & 0 &  0 &  0 \\
 0   & 0 &  0 &  0 & 0  & 0 &  0 &  0 
\end{bmatrix},
\end{equation}
\text{and}
 \begin{equation}
C = \begin{bmatrix}
 0   & 0 &  0 &  0 & 0  & 0 &  1 &  0 \\
 0   & 0 &  0 &  0 & 0  & 0 &  0 &  0 \\
 0   & 0 &  0 &  0 & 0  & 0 &  0 &  0 \\
 c^2 \, (\zeta_1 + \zeta_2 - \zeta_3)  & 0 &  0  & 0 & 0 & 0  & 0 & \zeta_1 \, \zeta_2 \\
 0   & 0 &  0 &  0 & 0  & 0 &  0 &  0 \\
 0   & 0 &  0 &  0 & 0  & 0 &  0 &  0 \\
 c^2 & 0 &  0 &  0 & 0  & 0 &  0 &  0 \\
 0   & 0 &  0 &  0 & 0  & 0 &  0 &  0 
\end{bmatrix}.
\end{equation}
By using a symbolic algebra program, we find that
the eigenvalues $\lambda$ of $P(i k)$ for (\ref{pml_sys_3d}) are 
\begin{equation}
\lambda \left( P \left(i k \right) \right) = \pm \, i \, c \, (k_1^2 + k_2^2 +  k_3^2  )^{1\slash 2}.
\end{equation}
Thus,
\begin{equation}
 \quad  \Re \{\lambda \left( P \left(i  k \right) \right) \} = 0, 
\end{equation}
while the corresponding eigenvectors are also complete, 
if at least two $\zeta_j = 0, \, j = 1, 2,3$; else, they are not complete.
Therefore, since $D$ is a diagonal matrix with negative entries for 
$\zeta_1, \, \zeta_2, \zeta_3 \geq 0$, we conclude that (\ref{pml-eq_3d}) is strongly stable, 
if at least two $\zeta_j = 0,$ and weakly stable, otherwise.
\section{Finite difference discretization}
Here we show how to discretize (\ref{pml-eq_3d}) with standard second-order finite differences
on a uniform mesh at grid points
$x_{i,i_{\ell}} = x_{1,0} + i_{\ell} \, \Delta x_{i}$, with $i=1,2,3$
and  $i_{\ell} = 0,1,\ldots, M_{\ell}$. For the time discretization we use a constant step
size $\Delta t$ and denote the time levels by $t_n = t_0 + n \Delta t$, $n = 0,1,\ldots, N$.
Inside the absorbing layer, we further introduce a space-time staggered grid at locations
$x_{i,i_{\ell}+\frac{1}{2}} = x_{i,0} + \left(i_{\ell} +\frac{1}{2}\right) \Delta x_{i}$, $i=1,2,3$ and times
$t_{n+\frac{1}{2}} = t_{0} + \left(n+\frac{1}{2}\right) \Delta t$. Then the numerical solution $u_{i ,j,k}^n$,
which approximates $u$ at grid point $(x_{1,i}, x_{2,j}, x_{3,k})$ and time $t_n$, satisfies
\begin{equation*}
\begin{split}
\, & \frac{u_{i,j,k}^{n+1} - 2 u_{i ,j,k}^n + u_{i,j,k}^{n-1}}{\Delta t \, ^2} 
+ \left(\zeta_{1 \, i} + \zeta_{2 \, j}  + \zeta_{3 \, k}  \right) \,\frac{u_{i,j,k}^{n+1} - u_{i,j,k}^{n-1}}{2 \Delta t } 
+ \left(\zeta_{1 \, i} \zeta_{2 \, j} + \zeta_{2 \, j} \zeta_{3 \, k} + \zeta_{3 \, k} \zeta_{1 \, i}   \right) u_{i ,j,k}^n \\
=  & \quad \frac{c^2_{i+\frac{1}{2},j,k}u^n_{i+1,j,k} -(c^2_{i+\frac{1}{2},j,k} + c^2_{i-\frac{1}{2},j,k} ) \, u^n_{i,j,k} + c^2_{i-\frac{1}{2},j,k}u^n_{i-1,j,k}}{\Delta x_1 \, ^2} \\
\, & + \frac{c^2_{i,j+\frac{1}{2},k}u^n_{i,j+1,k} -(c^2_{i,j+\frac{1}{2},k} + c^2_{i,j-\frac{1}{2},k} ) \, u^n_{i,j,k} + c^2_{i,j-\frac{1}{2},k}u^n_{i,j-1,k}}{\Delta x_2 \, ^2} \\
\, &+ \frac{c^2_{i,j,k+\frac{1}{2}}u^n_{i,j,k+1} -(c^2_{i,j,k+\frac{1}{2}} + c^2_{i,j,k-\frac{1}{2}} ) \, u^n_{i,j,k} + c^2_{i,j,k-\frac{1}{2}}u^n_{i,j,k-1}}{\Delta x_3 \, ^2}  \\
\, &+ \frac{\tilde{\phi}^n_{1 \,i+\frac{1}{2},j,k} -  \tilde{\phi}^n_{1 \,i-\frac{1}{2},j,k}}{\Delta x_1} + \frac{\tilde{\phi}^n_{2 \,i,j+\frac{1}{2},k} -  \tilde{\phi}^n_{2 \,i,j-\frac{1}{2},k}}{\Delta x_2} + \frac{\tilde{\phi}^n_{3 \,i,j,k+\frac{1}{2}} -  \tilde{\phi}^n_{3 \,i,j,k-\frac{1}{2}}}{\Delta x_3} -   \zeta_{1 \, i} \ \zeta_{2 \, j} \ \zeta_{3 \, k} \, \frac{\psi_{i,j,k}^{n+ \frac{1}{2}} +\psi_{i,j,k}^{n- \frac{1}{2}}  }{2},
\end{split}
\end{equation*}
where the cell averages of the auxiliary functions $\phi_{1}$, $\phi_{2}$ and $\phi_{3}$ are defined as
\begin{equation*}
\begin{split}
\tilde{\phi}^n_{1 \,i + \frac{1}{2},j,k} &= \frac{1}{4}\left(\phi_{1 \, i+\frac{1}{2},j - \frac{1}{2}, k - \frac{1}{2}}^n + \phi_{1 \, i+\frac{1}{2},j - \frac{1}{2}, k + \frac{1}{2}}^n + \phi_{1 \, i+\frac{1}{2},j + \frac{1}{2}, k - \frac{1}{2}}^n  + \phi_{1 \, i+\frac{1}{2},j + \frac{1}{2}, k + \frac{1}{2}}^n \right), \\
\tilde{\phi}^n_{2 \,i ,j+ \frac{1}{2},k} &= \frac{1}{4}\left(\phi_{2 \, i-\frac{1}{2},j + \frac{1}{2}, k - \frac{1}{2}}^n + \phi_{2 \, i-\frac{1}{2},j + \frac{1}{2}, k + \frac{1}{2}}^n + \phi_{2 \, i+\frac{1}{2},j + \frac{1}{2}, k - \frac{1}{2}}^n  + \phi_{2 \, i+\frac{1}{2},j + \frac{1}{2}, k + \frac{1}{2}}^n \right),  \\
\tilde{\phi}^n_{3 \,i ,j,k+ \frac{1}{2}} &= \frac{1}{4}\left(\phi_{3 \, i-\frac{1}{2},j - \frac{1}{2}, k + \frac{1}{2}}^n + \phi_{3 \, i-\frac{1}{2},j + \frac{1}{2}, k + \frac{1}{2}}^n + \phi_{3 \, i+\frac{1}{2},j - \frac{1}{2}, k + \frac{1}{2}}^n  + \phi_{3 \, i+\frac{1}{2},j + \frac{1}{2}, k + \frac{1}{2}}^n \right).
\end{split}
\end{equation*}
Concurrently with the above discretized wave equation, we also advance the (scalar) auxiliary variables 
$\psi$, $\phi_j$,  $j=1,2,3$ inside the absorbing layer by using standard 
finite differences. For $\psi$ we use
\begin{equation*}
\frac{\psi_{i,j,k}^{n+ \frac{1}{2}} - \psi_{i,j,k}^{n- \frac{1}{2}}}{\Delta t} =  u_{i,j,k}^n ,
\end{equation*}
whereas for $\phi_1$ we use
\begin{equation*}
\begin{split}
\,&\frac{\phi_{1 \, i+ \frac{1}{2}, j+ \frac{1}{2}, k+ \frac{1}{2}}^{n+1} - \phi_{1 \, i+ \frac{1}{2}, j+ \frac{1}{2}, k+ \frac{1}{2}}^n}{\Delta t}  \\
= &  - \zeta_{1 \, i+ \frac{1}{2}} \frac{\phi_{1 \, i+ \frac{1}{2}, j+ \frac{1}{2}, k+ \frac{1}{2}}^{n+1} + \phi_{1 \, i+ \frac{1}{2}, j+ \frac{1}{2}, k+ \frac{1}{2}}^n}{2} + \left(\zeta_{2 \, j+ \frac{1}{2}} + \zeta_{3 \, k+ \frac{1}{2}} - \zeta_{1 \, i+ \frac{1}{2}} \right) D_{x_1}^h u_{i+ \frac{1}{2}, j+ \frac{1}{2}, k+ \frac{1}{2}}^{n+ \frac{1}{2}} \\
\quad   &   + \zeta_{2 \, j+ \frac{1}{2}} \, \zeta_{3 \, k+ \frac{1}{2}}  D_{x_1}^h \psi_{1 \, i+ \frac{1}{2}, j+ \frac{1}{2}, k+ \frac{1}{2}}^{n+ \frac{1}{2}} \quad,
\end{split}
\end{equation*}
where 
\begin{equation*}
\begin{split}
 D_{x_1}^h u_{i+ \frac{1}{2}, j+ \frac{1}{2}, k+ \frac{1}{2}}^{n+ \frac{1}{2}} 
&= \frac{1}{2} \left(\frac{\tilde{u}^{n+1}_{i+1,j+ \frac{1}{2}, k+ \frac{1}{2}}- \tilde{u}^{n+1}_{i,j+ \frac{1}{2}, k+ \frac{1}{2}}    }{\Delta x_1} + \frac{\tilde{u}^{n}_{i+1,j+ \frac{1}{2}, k+ \frac{1}{2}}- \tilde{u}^{n}_{i,j+ \frac{1}{2}, k+ \frac{1}{2}}    }{\Delta x_1} \right), \\
\quad D_{x_1}^h \psi_{1 \, i+ \frac{1}{2}, j+ \frac{1}{2}, k+ \frac{1}{2}}^{n+ \frac{1}{2}}
&= \frac{\tilde{\psi}_{i+1,j+ \frac{1}{2},k+ \frac{1}{2}}^{n+\frac{1}{2}} -   \tilde{\psi}_{i, j+ \frac{1}{2},k+ \frac{1}{2} }^{n+\frac{1}{2}}}{ \Delta x_1}. 
\end{split}
\end{equation*}
Here, the cell averages of $u$ and $\psi$ are defined as
\begin{equation*}
\begin{split}
\tilde{u}^{n}_{i,j+ \frac{1}{2}, k+ \frac{1}{2}} &= \frac{1}{4}
\left(u_{i, j,k}^{n}+ u_{i, j,k+1}^{n} + u_{i, j+1,k}^{n} + u_{i,
  j+1,k+1}^{n} \right), \\ \tilde{\psi}_{i,j+
  \frac{1}{2},k+ \frac{1}{2}}^{n+\frac{1}{2}} &=
\frac{1}{4}\left(\psi_{i, j,k}^{n+\frac{1}{2}}+ \psi_{i,
  j,k+1}^{n+\frac{1}{2}} + \psi_{i, j+1,k}^{n+\frac{1}{2}} + \psi_{i,
  j+1,k+1}^{n+\frac{1}{2}} \right).
\end{split}
\end{equation*}
The finite difference approximations for $\phi_2$ and $\phi_3$ are analogous.
\section{Numerical experiments}
Here we present numerical experiments that illustrate the accuracy, versatility and long-time 
stability of our PML formulation discretized with standard finite differences
as in Section 4. In all cases we choose $\bar{\zeta_i} = 80$ in the damping profile,
which yields a relative reflection $R\approx 10^{-3}$ for the 
the typical values $c=1$ and $L_i = 0.1$. At the exterior boundary of the absorbing layer we
impose homogeneous Dirichlet boundary conditions.
\subsection{Point source in 2D}
First, we consider the wave equation (\ref{eq:wave1}) in two space
dimensions with constant speed of propagation $c=1$ and zero initial conditions, $u_0=v_0=0$.
The source term $f$ corresponds to a truncated first derivative of a Gaussian:
\begin{equation}
\label{eq:point_source}
f(x,y,t) = \delta (x) \, \delta (y) \,  h(t)
\end{equation}
with 
\begin{equation}
h(t) =  \frac{d}{dt} \, \big( \, \e^{-\pi^2  \left(f_0 \, t-1 \right)^2} \, \big), \qquad f_0 = 10 \,  \Hz.
\end{equation} 
The grid spacing is uniform in $x_1$ and $x_2$, with $\Delta x = 0.002$. \\
In Figure 2 we display snapshots of the numerical solutions 
at different times in $\Omega=[-0.5, \, 0.5]^2$, surrounded  by a PML of width $L=0.1$. We observe
how the circular wave propagates outward essentially without spurious reflection from the PML.
By time $t=1$ the wave has essentially left the computational domain.
To assess the error in the numerical solution, we compute a reference solution in a much larger
domain of size $[-5.5, 5.5]\times[5.5, 5.5]$, so that boundary effects are postponed to later times
inside $\Omega$. In Figure 3, the time evolution of the $L^2$--error is shown for different
values of the damping coefficient $\bar{\zeta_i}$. Until $t=8$ we observe a steady decrease 
of the error over seven orders of magnitude, regardless of the value of $\bar{\zeta_i}$,
which demonstrates the long-time stability of our method. Moreover,
our formulation appears robust with respect to the parameter value $\bar{\zeta_i}$.   
\begin{figure}[htbp]
\begin{center}
\begin{tabular}[t]{ccc}
\epsfig{file=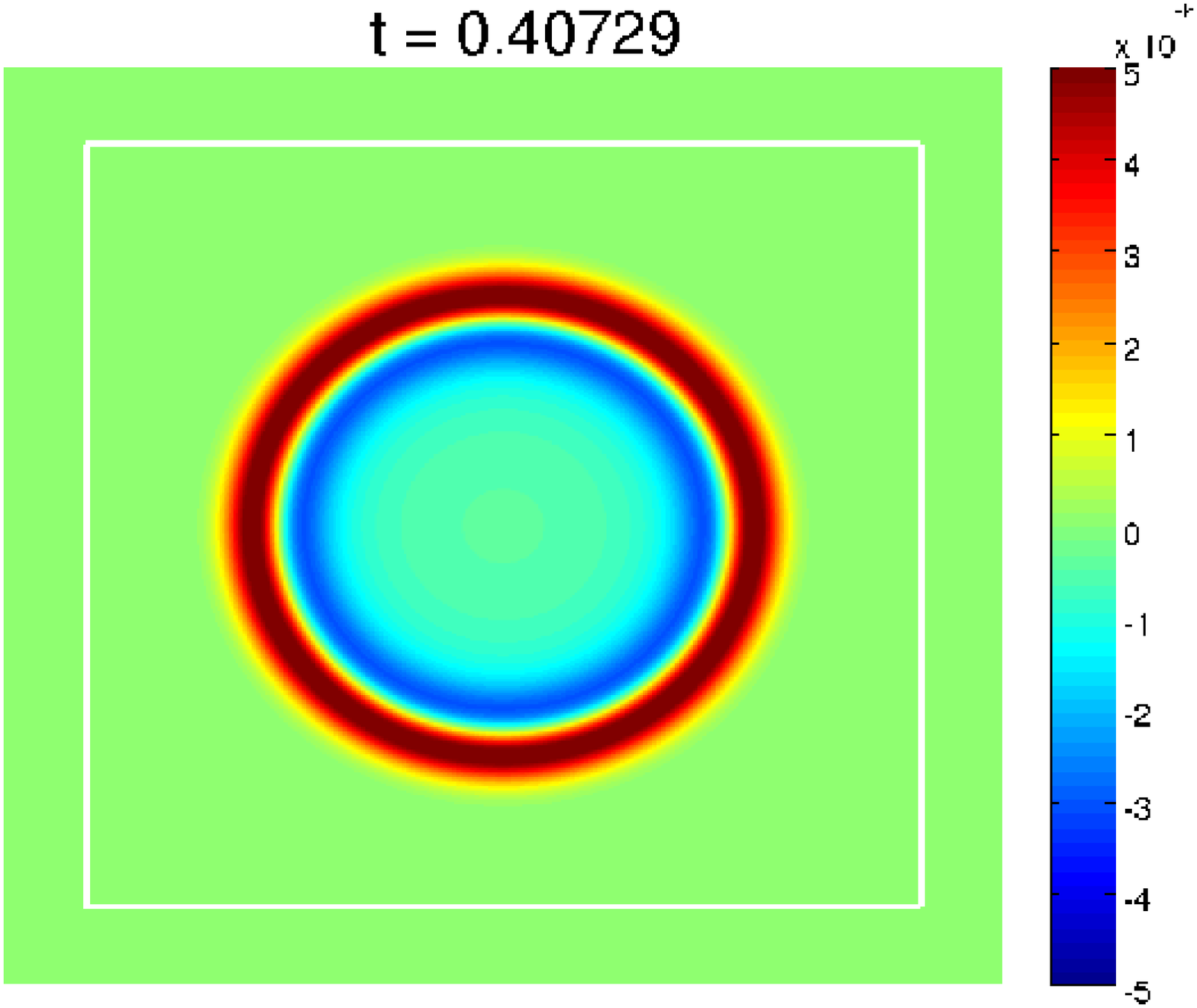, width=4.8cm} & \epsfig{file=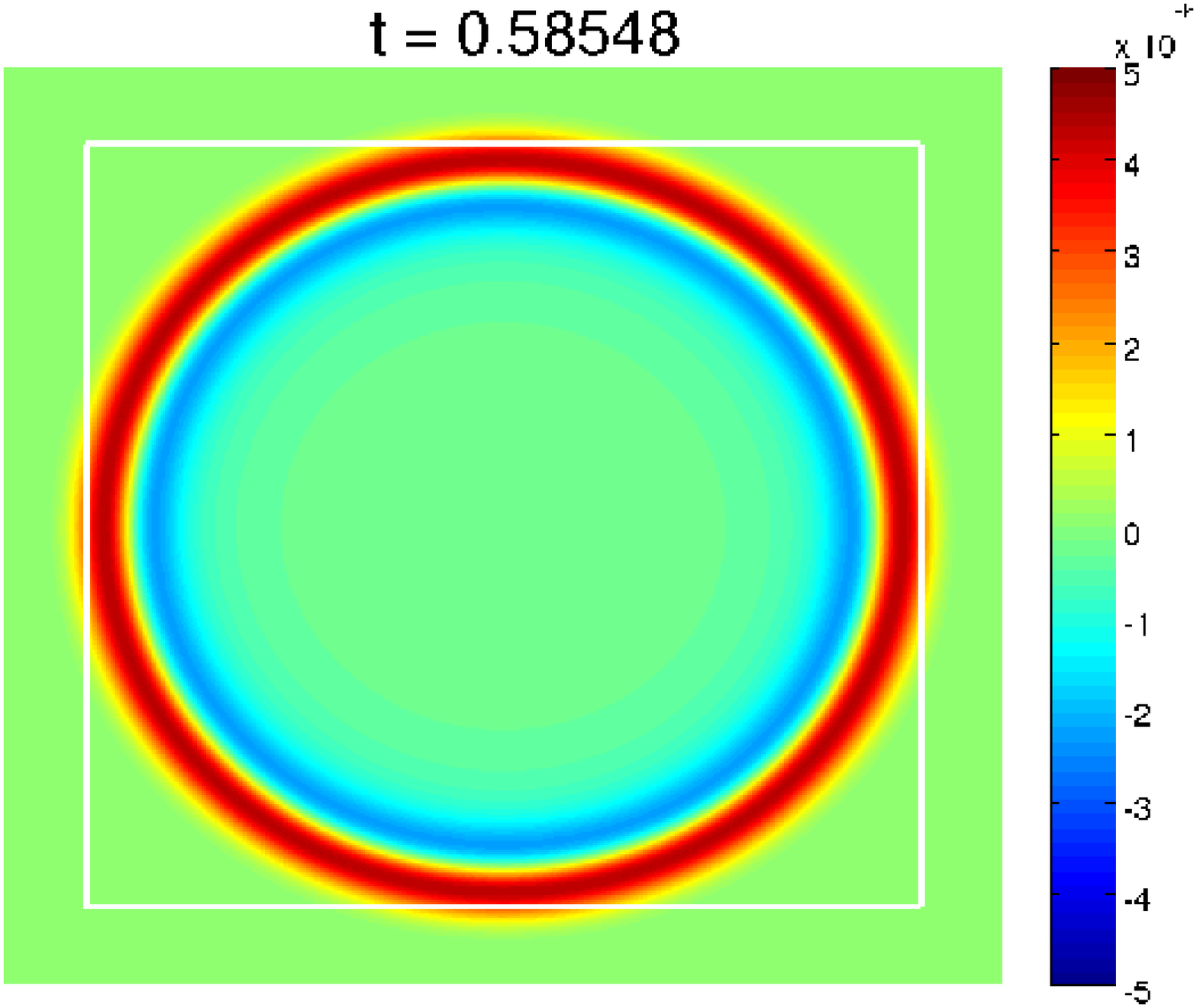, width=4.8cm} & \epsfig{file=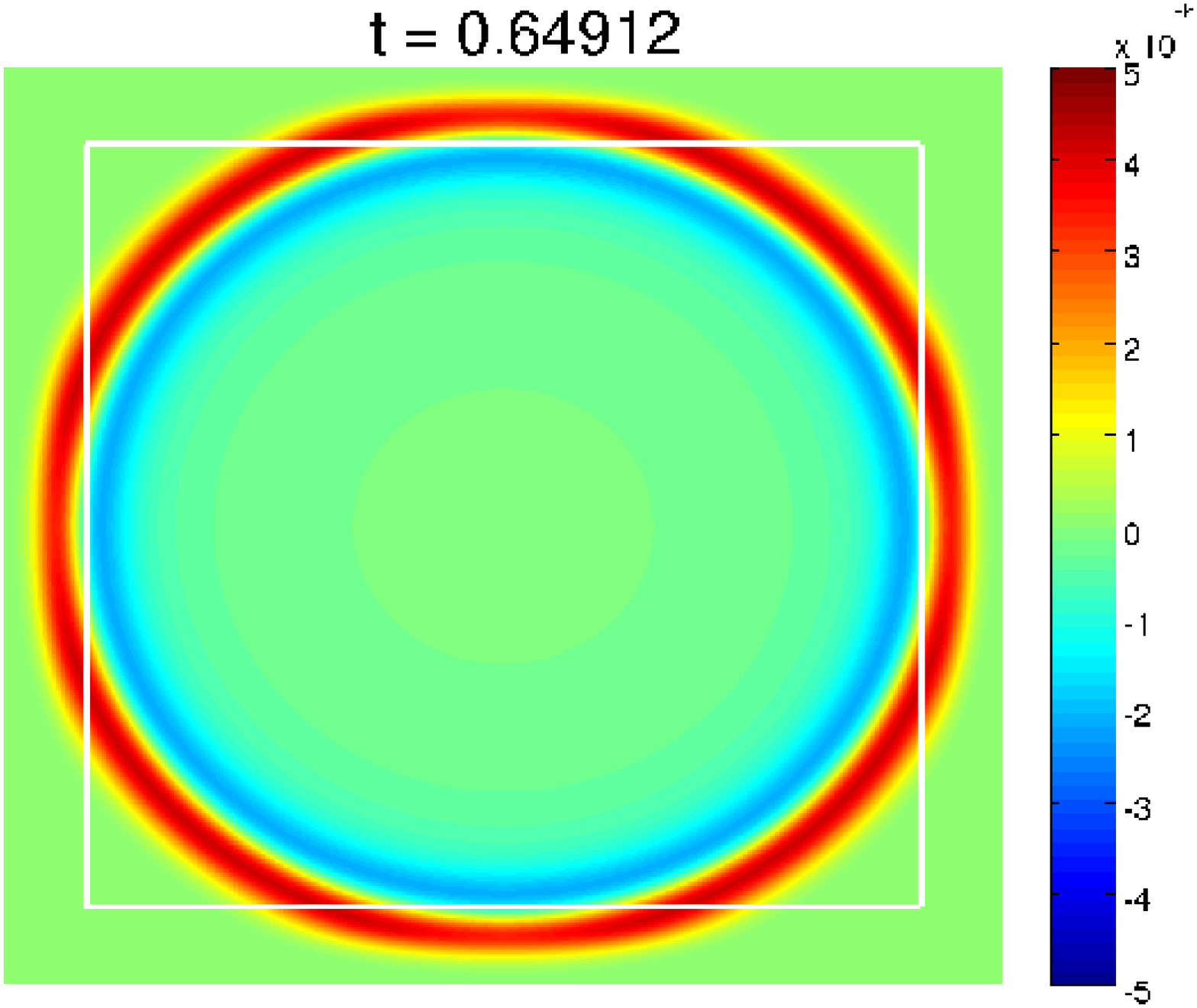, width=4.8cm}\\
\epsfig{file=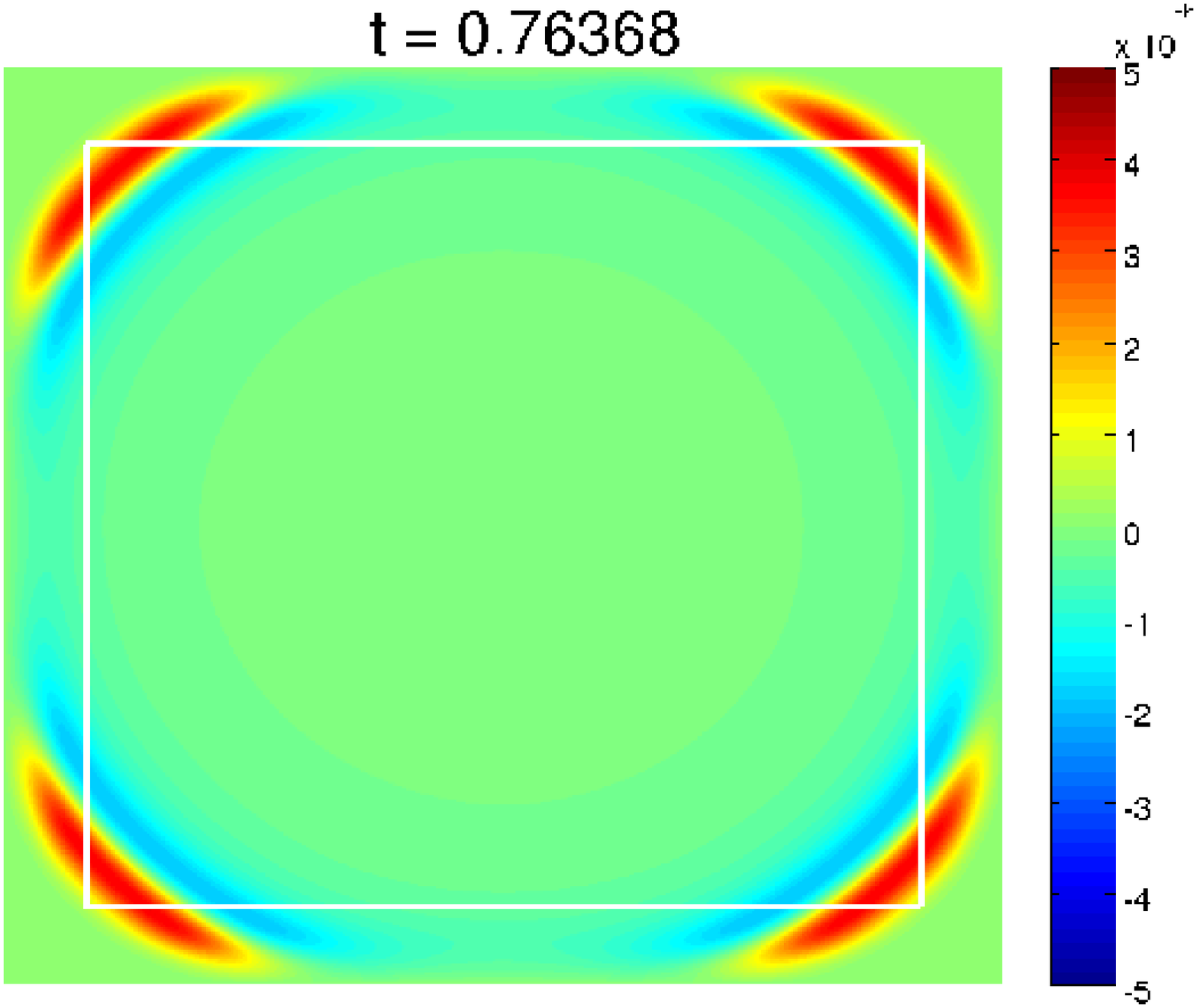, width=4.8cm} & \epsfig{file=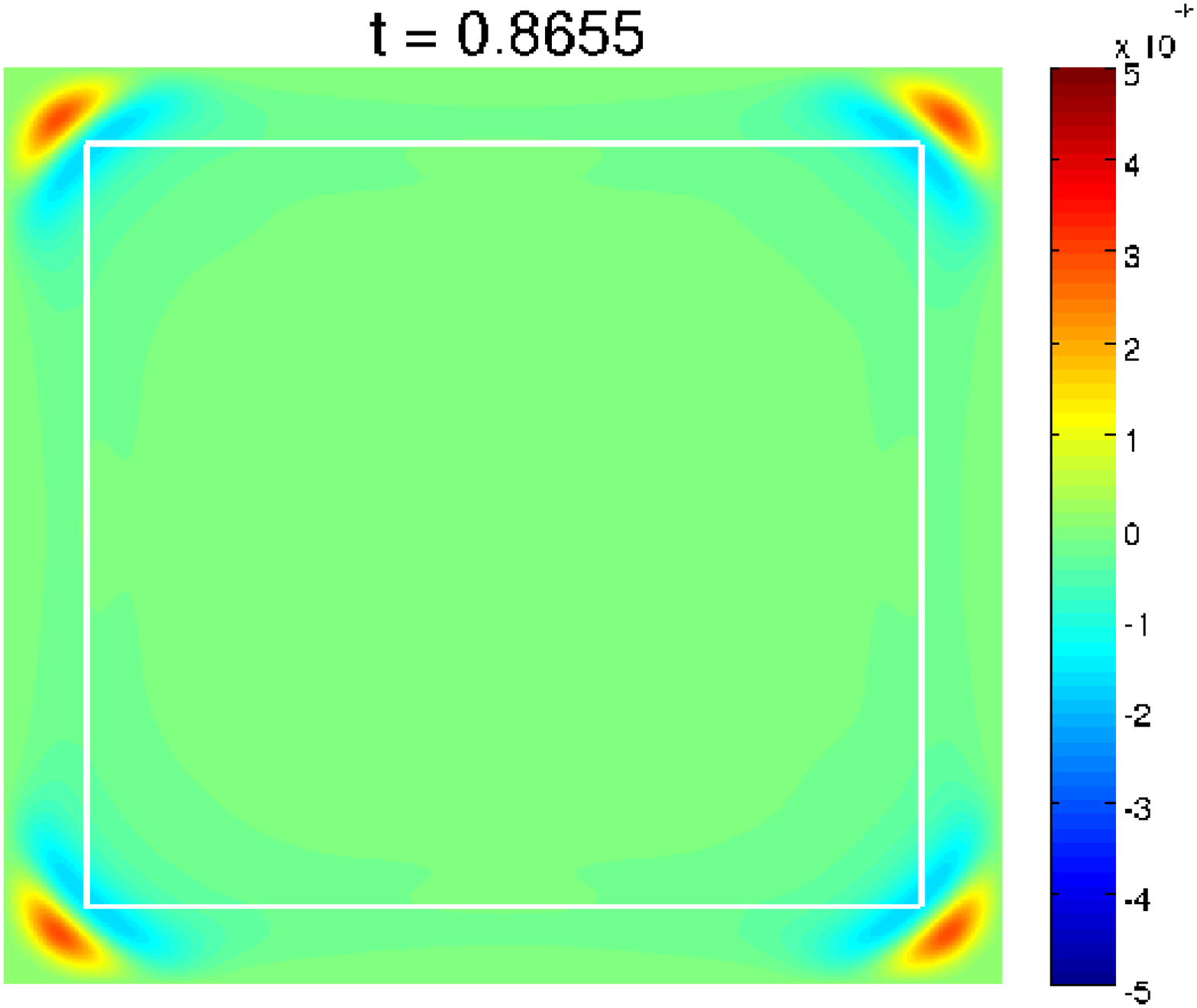, width=4.8cm} & \epsfig{file=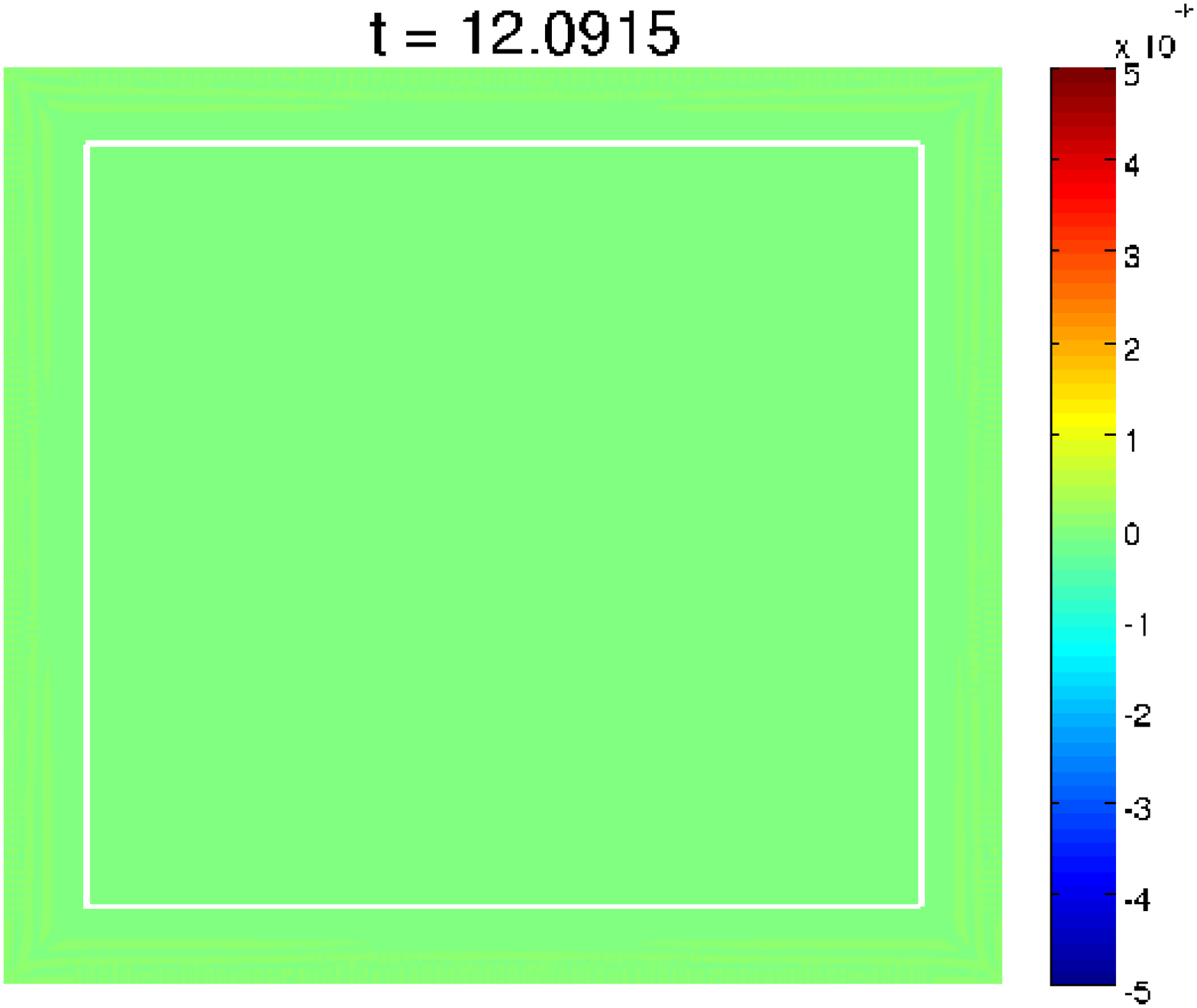, width=4.8cm}\\
\end{tabular}
\label{fig:gauss}
\end{center}
\caption{Point source in 2D: snapshots of the numerical solutions at different times in $\Omega=[-0.5, \, 0.5]^2$, surrounded  by a PML of width $L=0.1$.}
\end{figure}
\begin{figure}[htbp]
\begin{center}
\epsfig{file=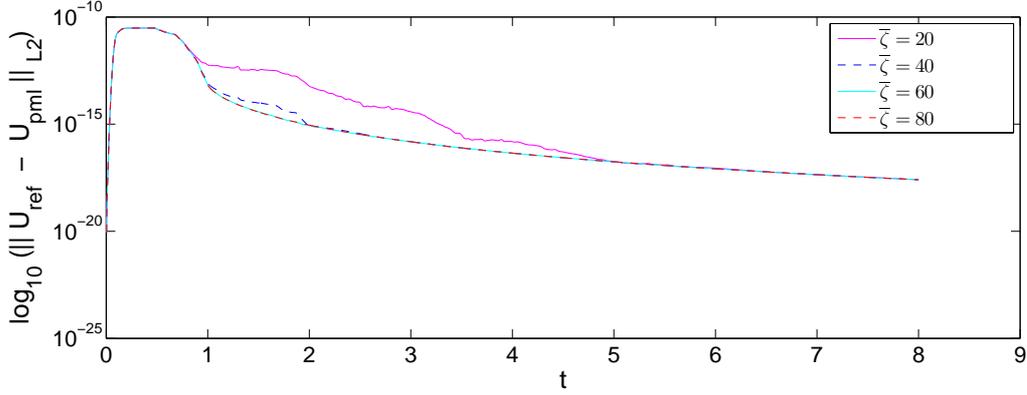, width=15.5cm} 
\end{center}
\caption{Point source in 2D: time evolution of the $L_2$--error for different damping coefficients 
$\bar{\zeta_i}$.}
\end{figure}
\subsection{Heterogeneous medium in 2D}
\begin{figure}[htbp]
\begin{center}
\epsfig{file=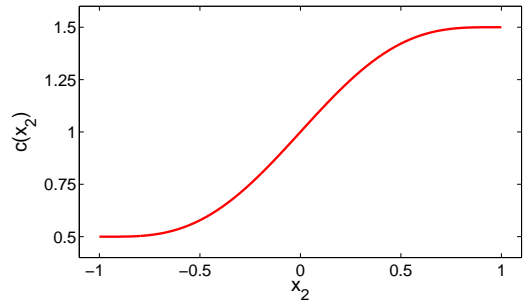, width=7.5cm}  
\end{center}
\caption{Heterogeneous medium in 2D: varying wave speed $c$ given by (\ref{eq:varying_c}).}
\label{fig:varying_c}
\end{figure}

\begin{figure}[htbp]
\begin{center}
\begin{tabular}[t]{ccc}
\epsfig{file=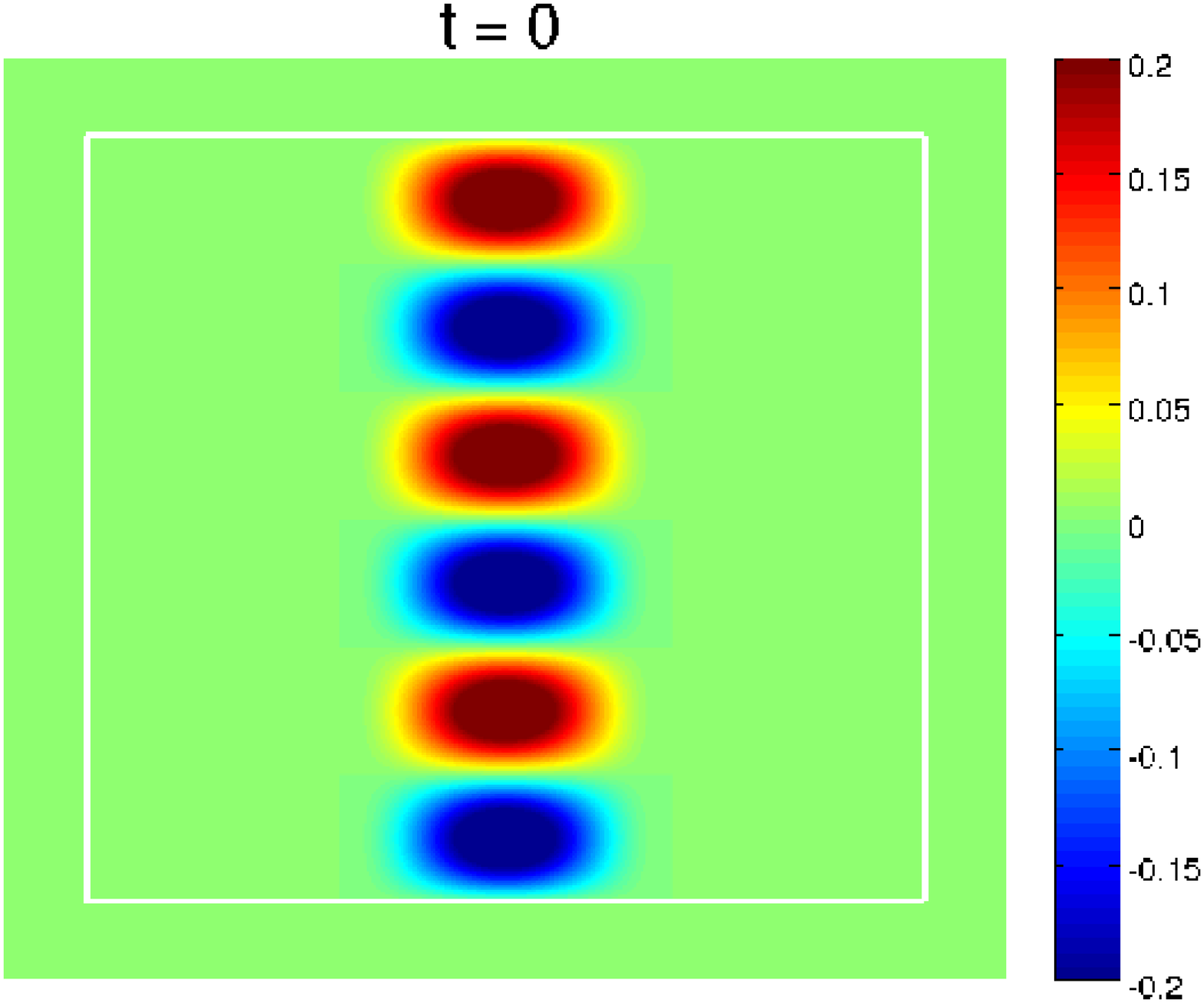, width=4.8cm} & \epsfig{file=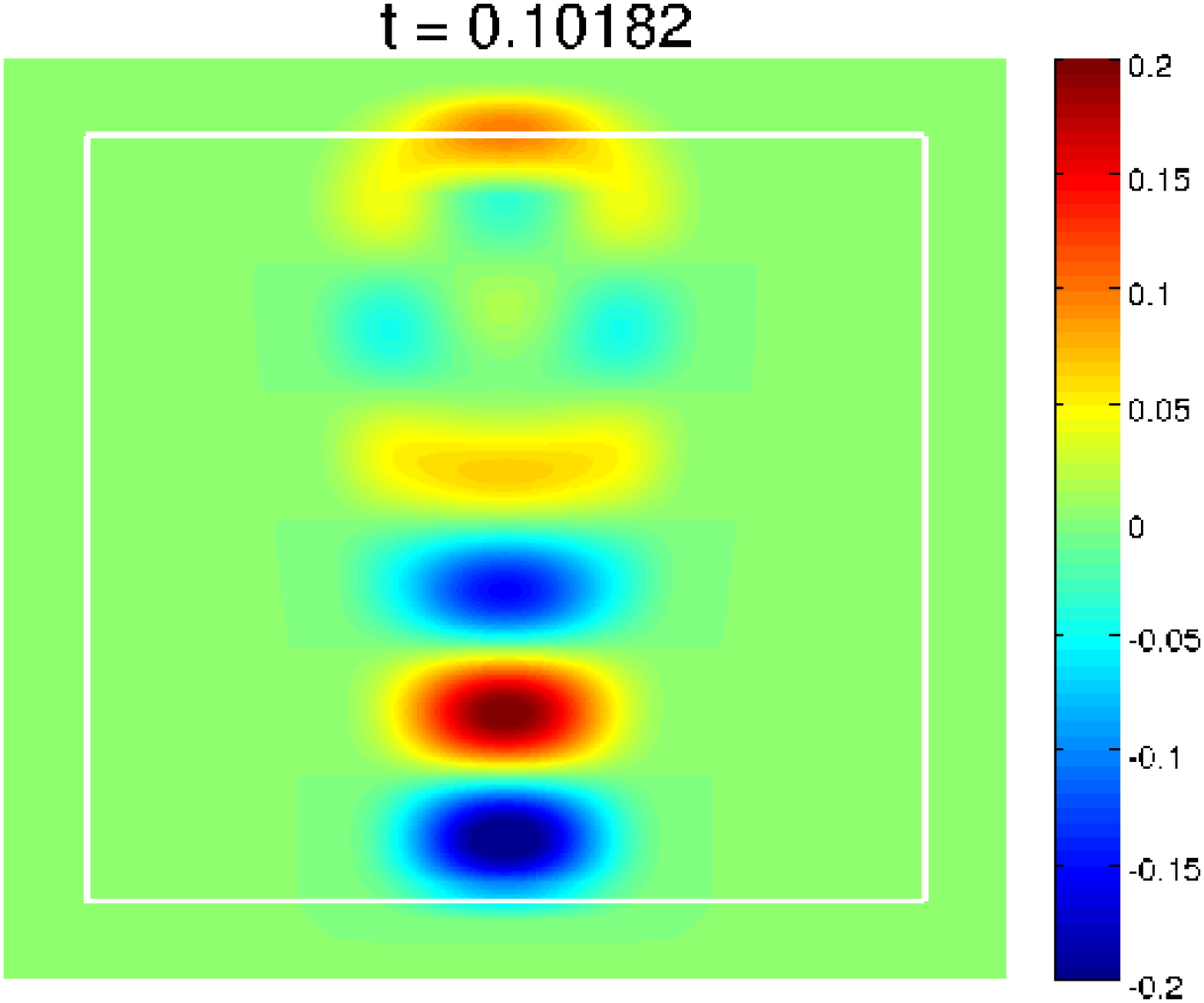, width=4.8cm} & \epsfig{file=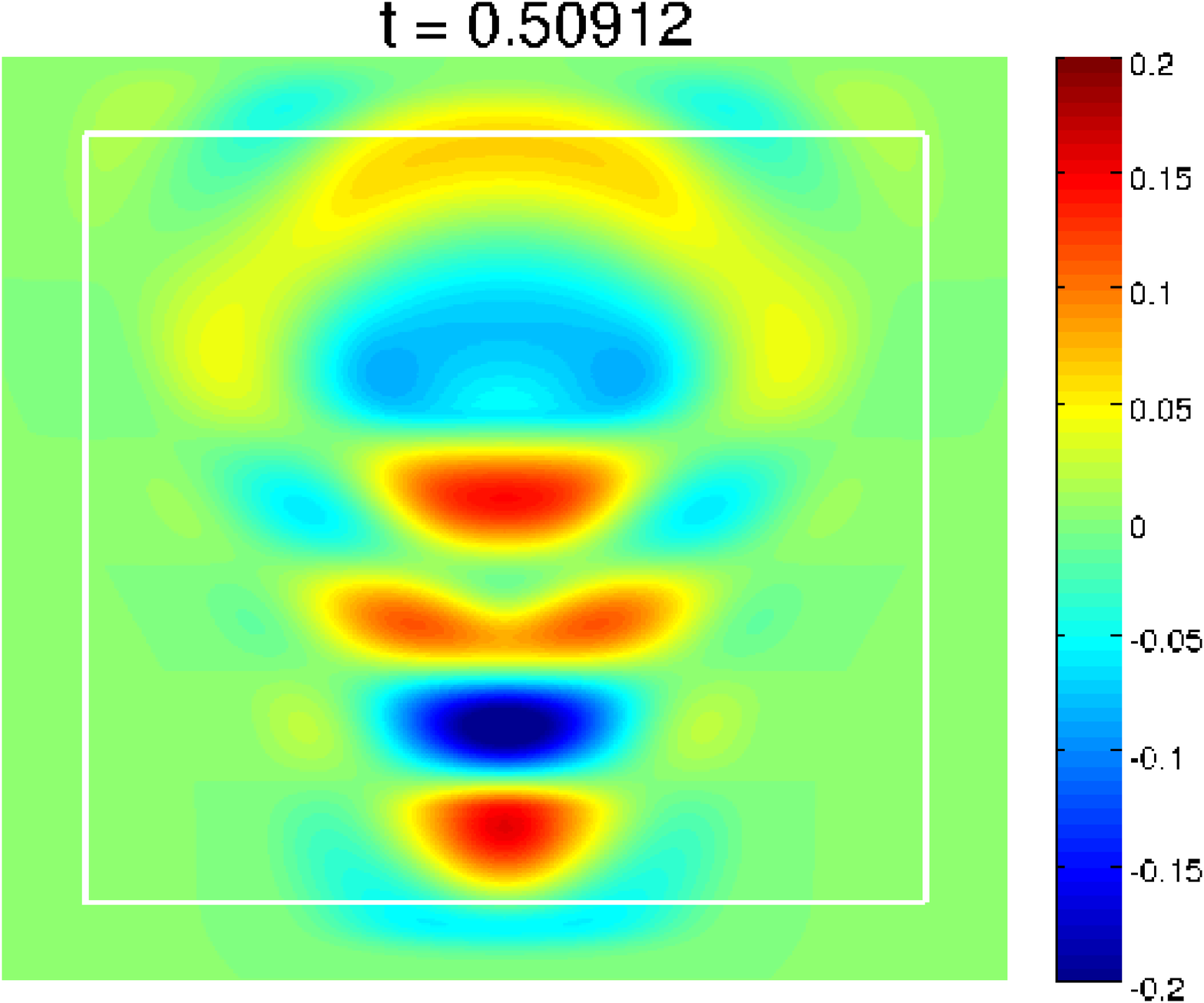, width=4.8cm} \\
\epsfig{file=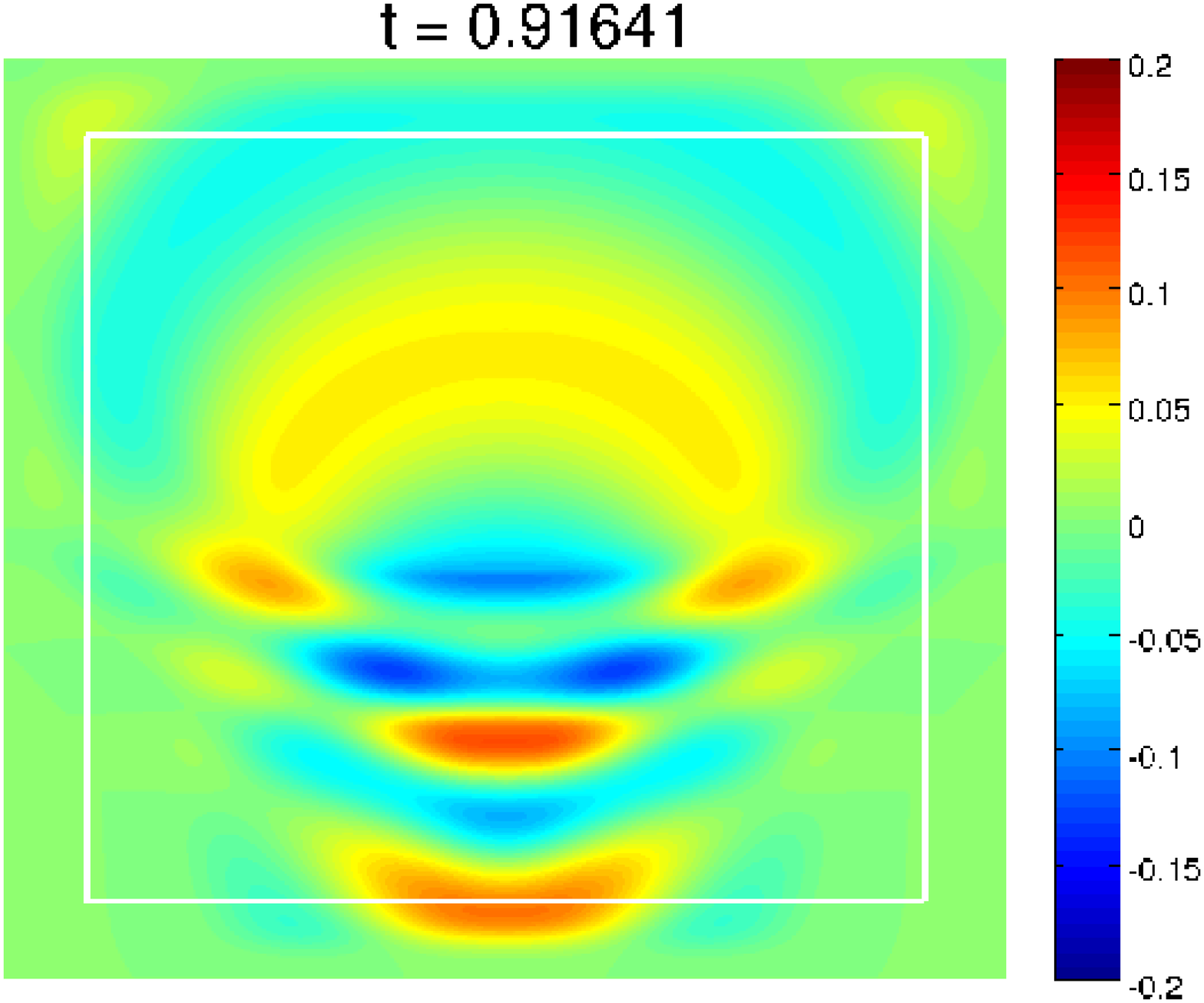, width=4.8cm} & \epsfig{file=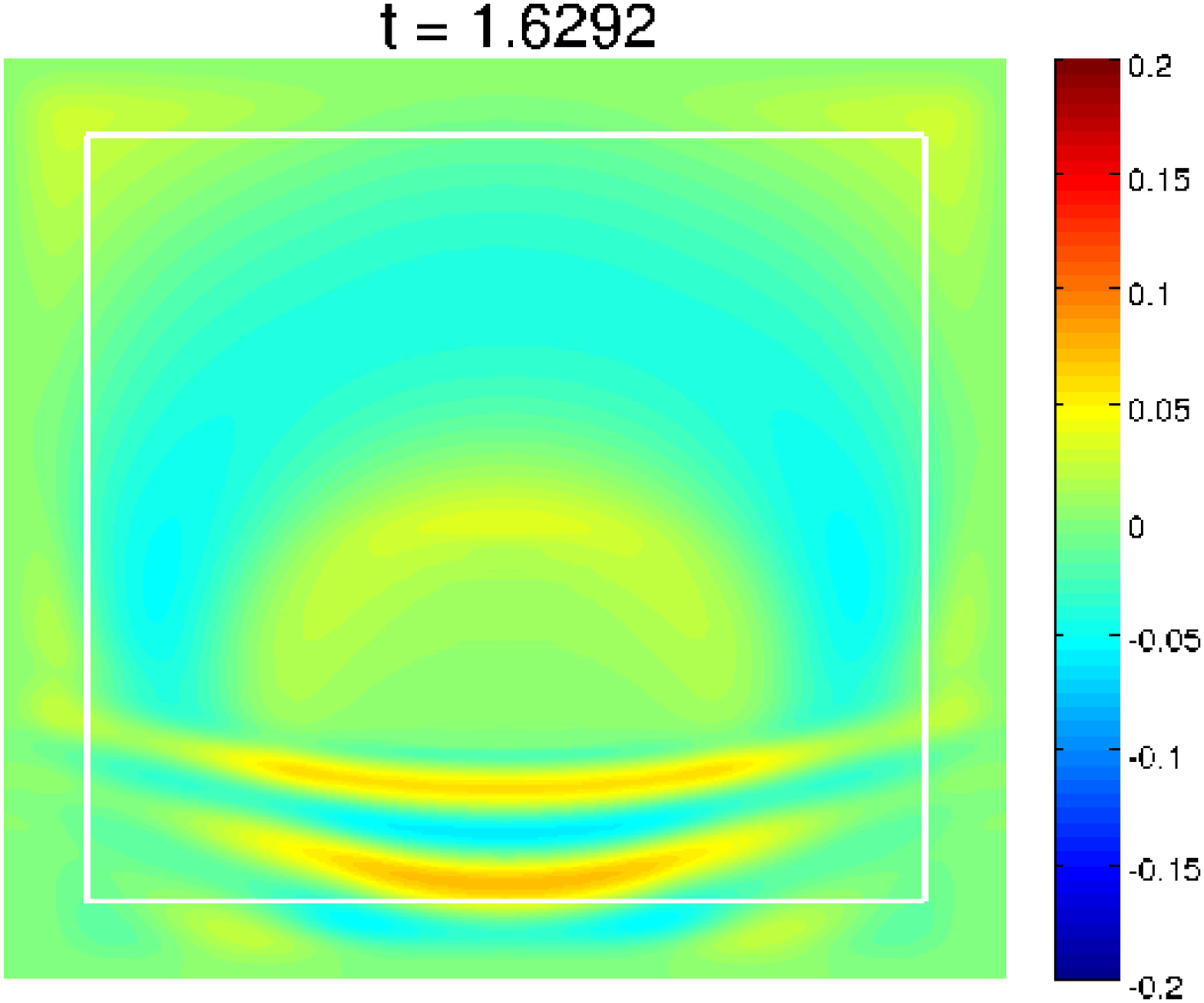, width=4.8cm} & \epsfig{file=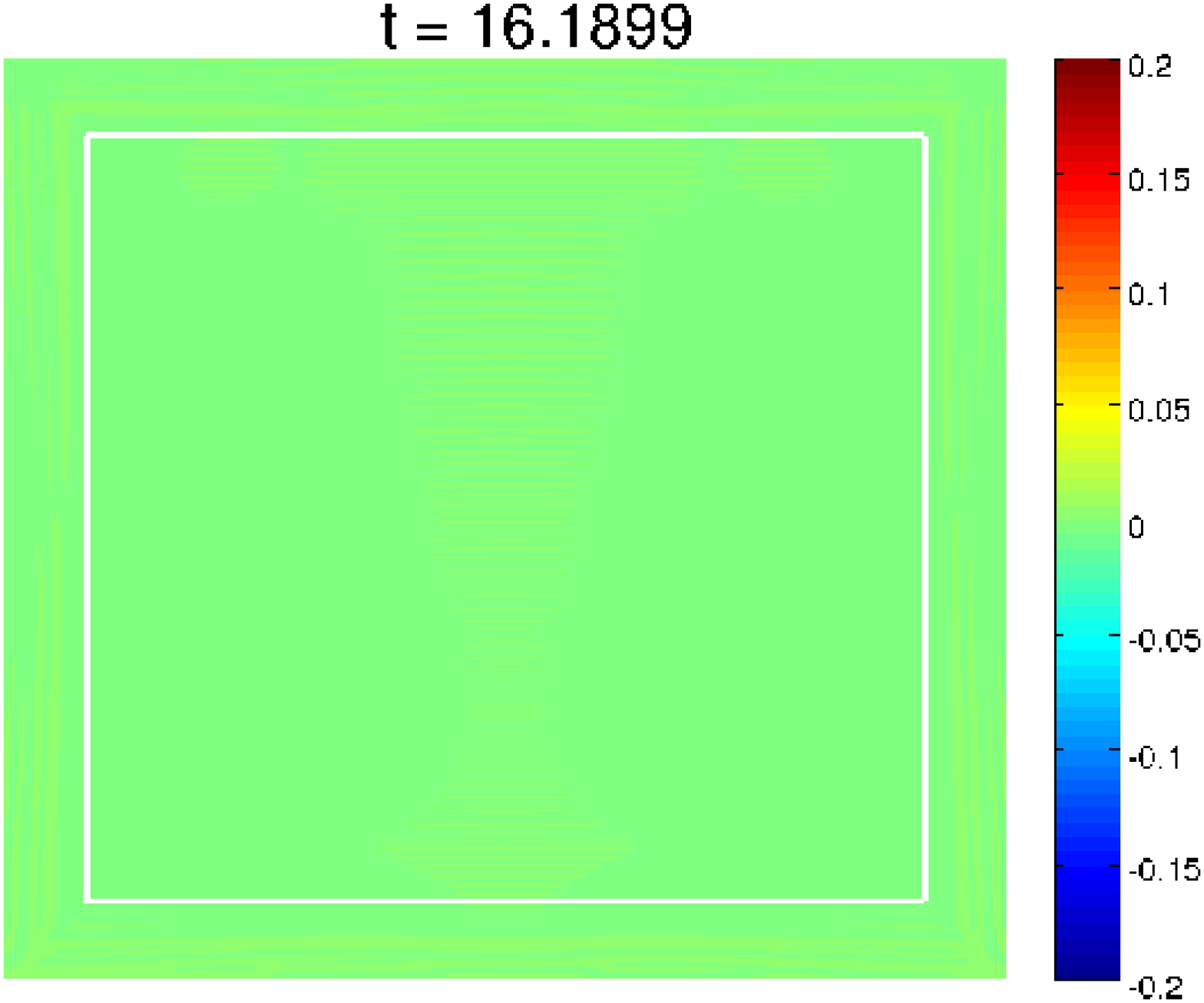, width=4.8cm} 
\end{tabular}
\end{center}
\caption{Heterogeous medium in 2D: snapshots of the numerical solution are shown at different times in 
$\Omega=[-1, \, 1]^2$, surrounded by a PML of width $L=0.2$.}
\label{fig:hetero}
\end{figure}

Next, to illustrate the versatility of our PML formulation, 
we consider the homogeneous wave equation (\ref{eq:wave1}) 
in a heterogeneous medium with varying wave speed $c=c(x_2)$, given by 
\begin{equation}
\label{eq:varying_c}
 c(x_1, x_2) = 
  \begin{cases}
     0.5,   &\text{if $ x_2< - b $} \\
     1 +    \frac{y}{2b} + \frac{1}{2\pi} \sin \left(\frac{\pi x_2}{b} \right), & \text{if $|x_2|<b $} \\
     1.5,   &\text{otherwise}. 
  \end{cases}
\end{equation}
We set $b=0.95$ which yields the vertical velocity profile shown in 
Figure \ref{fig:varying_c}.
The initial conditions are
\begin{equation}
u|_{t=0} = u_0(x_1,x_2) \quad \text{and} \quad   u_t|_{t=0} = 0 ,
\end{equation}
where 
\begin{equation*}
u_0(x_1,x_2) =
  \begin{cases}
     \big(4 \left(x_1+0.4\right) \left(0.4-x_1 \right) \big)^3 \sin(3\pi  x_2), &\text{if $-0.4<x_1<0.4, \, -1 < x_2 <1$} \\
      0,                                            &\text{otherwise.}
  \end{cases}
\end{equation*}
Here $\Omega$ is the square domain $[-1, \, 1] \times [-1, \, 1]$, surrounded by a PML of width $L=0.2$. 
The finite difference grid is uniform with grid spacing $\Delta x = 0.004$.
In Fig.\ref{fig:hetero}, we display snapshots of the solution at different times, where
again the last frame is purposely chosen at a much later time. In spite of the
varying wave speed and the glancing angle of incidence along the vertical artificial
boundaries, the waves are damped without spurious reflection. Even at much later
times we do not observe any instability in the numerical scheme.
\subsection{Point source in 3D}
Finally, we consider the wave equation (\ref{eq:wave1}) in three space dimensions
with zero initial conditions and the same point source $f$ as in  
(\ref{eq:point_source}). The grid spacing is uniform 
in $x_1$, $x_2$ and $x_3$ with $\Delta x = 0.006$.
In Figure 6, we display snapshots of the numerical solutions 
at different times in $\Omega=[-0.5, \, 0.5]^2$, surrounded  by a PML of width $L=0.1$. We observe
how the spherical wave propagates outward essentially without spurious reflection from the PML.
By time $t=1$ the wave has essentially left the computational domain. Again we observe
no instabilities in the numerical solution even at much later times.
\begin{figure}[htbp]
\begin{center}
\begin{tabular}[t]{cc}
\epsfig{file=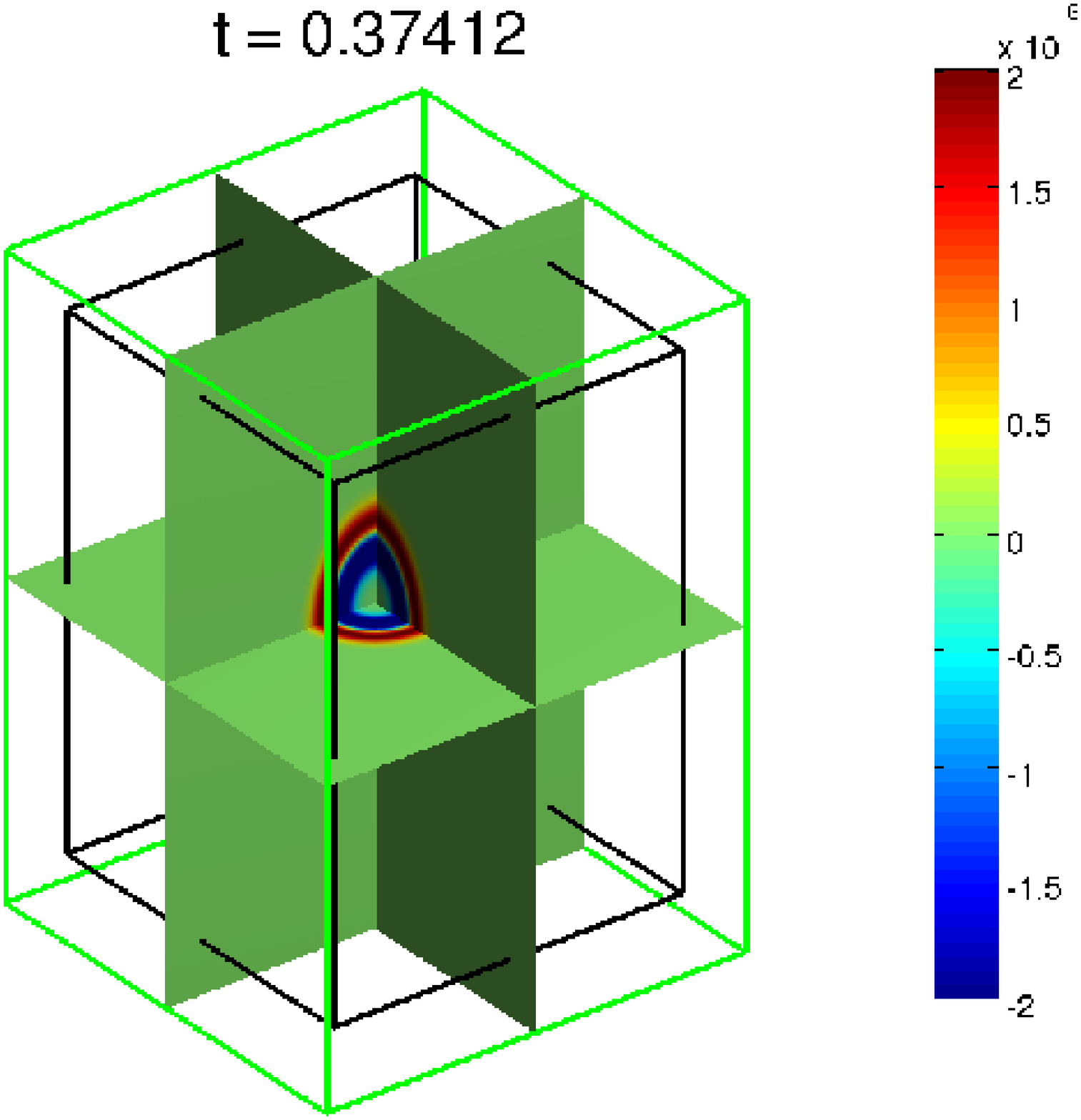, width=8.7cm} & \epsfig{file=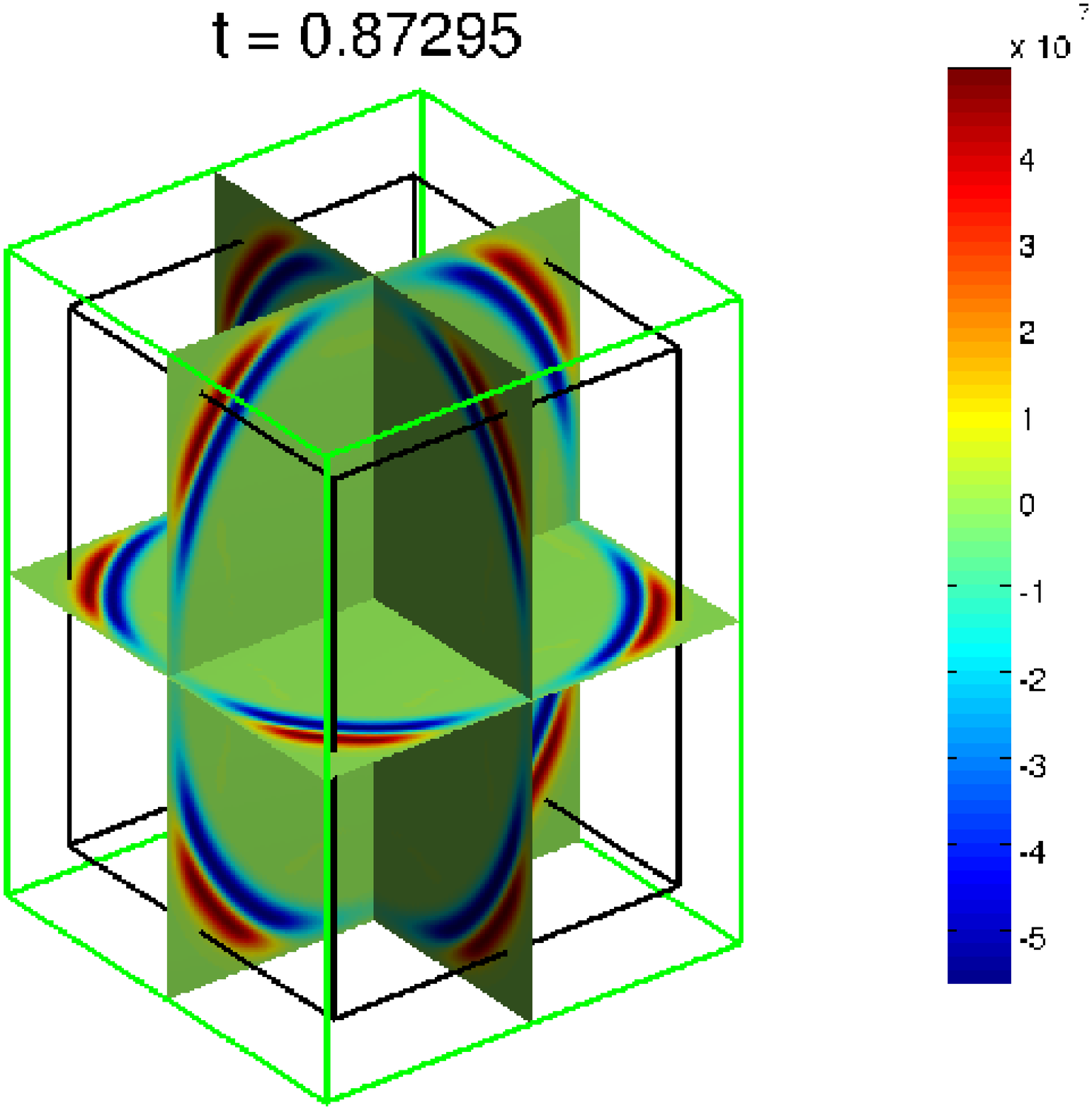, width=8.7cm} \\
\epsfig{file=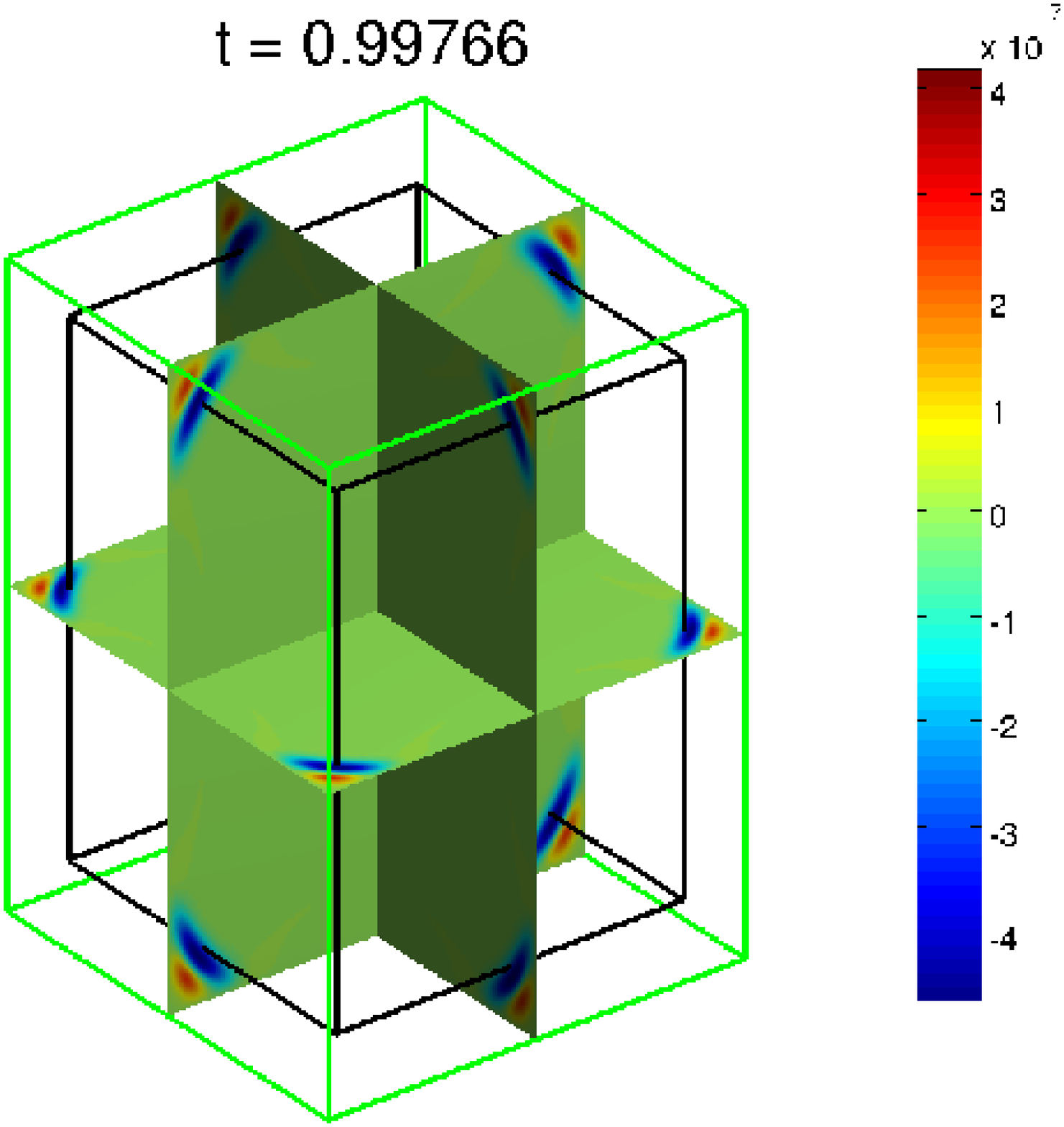, width=8.7cm} & \epsfig{file=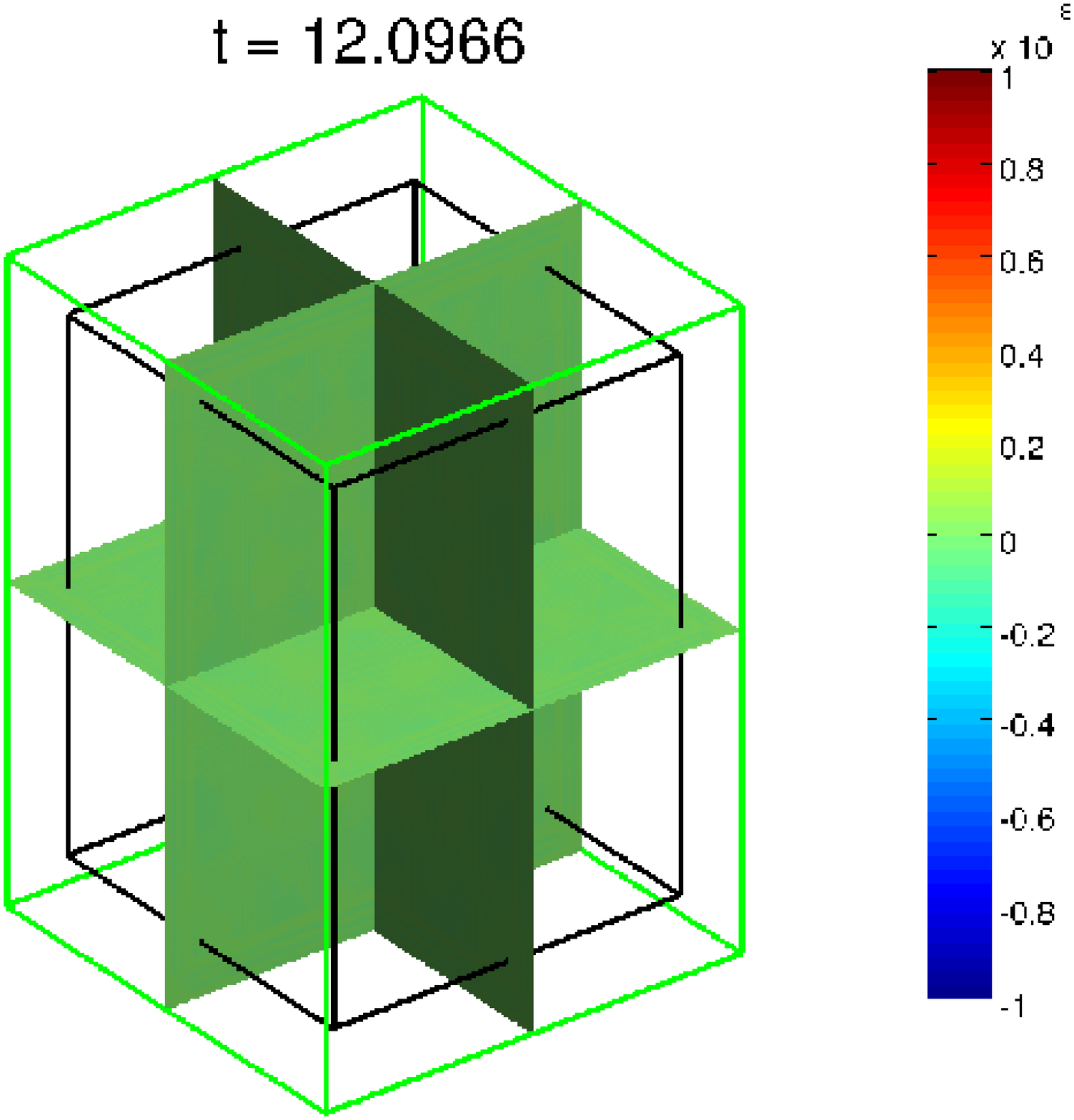, width=8.7cm} 
\end{tabular}
\label{fig:gauss3d}
\end{center}
\caption{Point source in 3D: snapshots of the numerical solution are shown at different times in 
$\Omega=[-0.5, \, 0.5]^3$, surrounded  by a PML of width $L=0.1$.}
\end{figure}
\section{Concluding remarks}
We have presented a PML formulation for the wave equation in its standard second-order form. 
It distinguishes itself from known formulations by its simplicity and the small number of auxiliary 
variables needed inside the absorbing layer. We have proved that the continuous
Cauchy problem with the unbounded PML is stable and well-posed. Our numerical results in two
and in three space dimensions with standard finite differences 
illustrate the accuracy, versatility and long-time stability of our PML formulation. \\
Because it involves no high space or time derivatives, our PML formulation easily fits continuous
or discontinuous Galerkin formulation for use with finite element methods \cite{dg,gss}. It also immediately
generalizes to Maxwell's equations in second-order form. Current work involves
the extension to second-order wave equations in complex elastic and poro-elastic media,
and will be reported elsewhere in the near future.

\end{document}